\documentstyle[11pt,amssymb,amsmath]{article}

\newcommand{\be}{\begin{equation}}
\newcommand{\ee}{\end{equation}}
\newcommand{\bea}{\begin{eqnarray}}
\newcommand{\eea}{\end{eqnarray}}
\newtheorem{theorem}{Theorem}

\newtheorem{conjecture}{Conjecture}
\newtheorem{corollary}{Corollary}

\newtheorem{lemma}{Lemma}

\def\1#1{^{(#1)}}

\begin{document}
\title{Arnold's Conjectures on Weak Asymptotics\\and\\Statistics of
Numerical Semigroups ${\sf S}\left(d_1,d_2,d_3\right)$}
\author{Leonid G. Fel\\
\\Department of Civil Engineering, 
Technion, Haifa 3200, Israel\\
\vspace{-.3cm}  
\\{\sl e-mail: lfel@techunix.technion.ac.il}}
\date{\today}
\maketitle

\def\be{\begin{equation}}
\def\ee{\end{equation}}  
\def\bea{\begin{eqnarray}}   
\def\eea{\end{eqnarray}}
\def\p{\prime}
\begin{center}
\end{center}

\begin{abstract}
Three conjectures $\#$1999--8, $\#$1999--9 and $\#$1999--10 which were posed 
by V. Arnold \cite{arn00} and devoted to the statistics of the numerical 
semigroups are refuted for the case of semigroups generated by three positive 
integers $d_1,d_2,d_3$ with $\gcd(d_1,d_2,d_3)=1$. Weak asymptotics of 
conductor $C(d_1,d_2,d_3)$ of numerical semigroup and fraction $p(d_1,d_2,d_3)$ 
of a segment [$0;C(d_1,d_2,d_3)-1$] occupied by semigroup are found.
\end{abstract}

\newpage
\tableofcontents
\newpage
\section{Introduction}\label{intr}
Some years ago V. Arnold has posed three conjectures \cite{arn99}, \cite{arn00},
\cite{arn04} devoted to statistics of numerical semigroups generated by $m$ 
positive integers $d_1,\ldots,d_m$ for $m\geq 3$. This statistics is concerned 
with the so--called {\em weak asymptotics} for the numbers of the integer 
points in the domains and on the surface in ${\mathbb R}^m$. The conjectures 
are enumerated in \cite{arn00} as $\#$1999--8, $\#$1999--9 and $\#$1999--10 and 
are intimately related to the Frobenius problem for the numerical semigroups 
where a progress was achieved recently \cite{fel04} in the case $m=3$. This is 
the first nontrivial case where a set of numerical semigroups is separated into 
symmetric and non--symmetric semigroups with rather different homological 
properties of their associated polynomial rings \cite{frob87}, \cite{brown91}, 
\cite{brown92}. Based on these properties we refute the Arnold's conjectures 
for semigroups generated by three elements.

The paper is organized in six Sections. In Section \ref{basic0} we recall the 
main facts about numerical semigroups generated by three elements and their 
associated polynomial rings. Following \cite{arn04}, in Section \ref{stat1} we 
define a weak asymptotic of numerical functions on semigroups at the typical 
large vectors. In Section \ref{numer1} we prove a technical Lemma \ref{lem2} 
on statistics of symmetric and non--symmetric semigroups generated by three 
elements which makes a basis to perform calculations in the following Sections. 
In Section \ref{arr1} we refute the conjectures of Arnold for semigroups 
generated by three elements. In Section \ref{rev1} we show that two weak 
asymptotics, for conductor $C$ of semigroup and fraction $p$ of a segment 
[$0;C-1$] occupied by semigroup, are not universal and depend on typical vectors
where an averaging is performed around. Based on results of \cite{fel04} we 
also improve the lower bound of $p$ which was obtained in Section \ref{arr1} 
with less powerful methods. 
\section{Algebra of numerical semigroups ${\sf S}\left(d_1,d_2,d_3\right)$}
\label{basic0}
Let ${\sf S}\left(d_1,d_2,d_3\right)\subset {\mathbb Z}_+$ be the additive 
numerical semigroup finitely generated by a minimal set of positive integers 
$\{d_1,d_2,d_3\}$ such that $d_1<d_2<d_3$ and $\gcd(d_1,d_2,d_3)=1$. It is 
classically known that $d_1\geq 3$ \cite{abhy67}. For short we denote the 
vector $(d_1,d_2,d_3)$ by ${\bf d}^3$. The least positive integer ($d_1$)
belonging to ${\sf S}\left({\bf d}^3\right)$ is called {\em the multiplicity}. 
The smallest integer $C\left({\bf d}^3\right)$ such that all integers $s,\;s
\geq C\left({\bf d}^3\right)$, belong to ${\sf S}\left({\bf d}^3\right)$ is 
called {\em the conductor} of ${\sf S}\left({\bf d}^3\right)$,
\begin{eqnarray}
C\left({\bf d}^3\right):=\min\left\{s\in {\sf S}\left({\bf d}^3\right)\;|\;
s+{\mathbb Z}_+\cup\{0\}\subset {\sf S}\left({\bf d}^3\right)\right\}\;.
\label{intro1}
\end{eqnarray}
The number $F\left({\bf d}^3\right)=C\left({\bf d}^3\right)-1$ is referred 
to as {\em the Frobenius number}. Denote by $\Delta\left({\bf d}^3\right)$ 
the complement of ${\sf S}\left({\bf d}^3\right)$ in ${\mathbb Z}_+$, i.e. 
$\Delta\left({\bf d}^3\right)={\mathbb Z}_+\setminus {\sf S}\left({\bf d}^3
\right)$. The cardinalities ($\#$) of the set $\Delta\left({\bf d}^3\right)$
and the set ${\sf S}\left({\bf d}^3\right)\cap $[0;$F\left({\bf d}^3\right)$]
are called {\em a number of gaps} $G\left({\bf d}^3\right)$, or {\em a genus}
of ${\sf S}\left({\bf d}^3\right)$, and {\em a number of nongaps}
$\widetilde{G}\left({\bf d}^3\right)$, respectively,
\begin{eqnarray}
&&G\left({\bf d}^3\right):=\#\left\{\Delta\left({\bf d}^3\right)\right\}\;,
\;\;\;\widetilde{G}\left({\bf d}^3\right):=\#\left\{{\sf S}\left({\bf d}^3
\right)\cap [0; F\left({\bf d}^3\right)]\right\}\;,\;\;\;\mbox{so that}
\label{intro2a}\\
&&G\left({\bf d}^3\right)+\widetilde{G}\left({\bf d}^3\right)=
C\left({\bf d}^3\right)\;.\label{intro2b}
\end{eqnarray}
Notice that two requirements, $\gcd(d_1,d_2,d_3)=1$ and $G\left({\bf d}^3
\right)<\infty$, are equivalent. 

The semigroup ${\sf S}\left({\bf d}^3\right)$ is called {\em symmetric} iff 
for any integer $s$ holds
\begin{eqnarray}
s\in {\sf S}\left({\bf d}^3\right)\;\;\;\Longleftrightarrow\;\;\;F\left({
\bf d}^3\right)-s\not\in{\sf S}\left({\bf d}^3\right)\;.\label{intro3}
\end{eqnarray}
Otherwise ${\sf S}\left({\bf d}^3\right)$ is called {\em non--symmetric}. The 
integers $G\left({\bf d}^3\right)$ and $C\left({\bf d}^3\right)$ are related 
as \cite{heku71},
\begin{eqnarray}
2G\left({\bf d}^3\right)=C\left({\bf d}^3\right)\;\;\mbox{if}\;\;{\sf S}
\left({\bf d}^3\right)\;\;\mbox{is symmetric semigroup, and}\;\;2G\left({\bf 
d}^3\right)>C\left({\bf d}^3\right)\;\;\mbox{otherwise}.\label{intro5}
\end{eqnarray}
Notice that ${\sf S}\left({\bf d}^2\right)$ is always symmetric semigroup
\cite{aper46}. 

Let $R={\sf k}\left[X_1,X_2,X_3\right]$ be a polynomial ring in 3 variables 
over a field ${\sf k}$ of characteristic 0 and
\begin{eqnarray}
\pi\;:\;\;\;{\sf k}\left[X_1,X_2,X_3\right]\longmapsto {\sf k}\left[z^{
d_1},z^{d_2},z^{d_3}\right]\nonumber
\end{eqnarray}
be the projection induced by $\pi\left(X_i\right)=z^{d_i}$. Denote ${\sf k}
\left[z^{d_1},z^{d_2},z^{d_m}\right]$ by ${\sf k}\left[{\sf S}\left({\bf d}^3
\right)\right]$. Then ${\sf k}\left[{\sf S}\left({\bf d}^3\right)\right]$ is
a graded subring of ${\sf k}\left[X_1,X_2,X_3\right]$ and has a presentation
as a R--module,
\begin{eqnarray}
{\sf k}\left[{\sf S}\left({\bf d}^3\right)\right]\cong {\sf k}\left[X_1,X_2,
X_3\right]/{\cal I}\left({\bf d}^3\right)\;.\nonumber
\end{eqnarray}
The prime ideal ${\cal I}\left({\bf d}^3\right)$ is the kernel of the map 
$\pi$ and it is minimally generated by a finite number of generators $P_k\left(
X_1,X_2,X_3\right)$ such that $\pi\left(P_k\right)=P_k\left(z^{d_1},z^{d_2},
z^{d_3}\right)=0$. 

The ring ${\sf k}\left[{\sf S}\left({\bf d}^3\right)\right]$ is a 1--dim 
Cohen--Macaulay ring \cite{kunz85}, and becomes Gorenstein ring iff ${\sf S}
\left({\bf d}^3\right)$ is symmetric \cite{kunz70}. Moreover, by \cite{serr63} 
the Gorenstein ring ${\sf k}\left[{\sf S}\left({\bf d}^3\right)\right]$ is a 
complete intersection. 
Denote by $t\left({\sf S}\left({\bf d}^3\right)\right)$ {\em a type} of the 
ring ${\sf k}\left[{\sf S}\left({\bf d}^3\right)\right]$ which in the case of 
numerical semigroup coincides with a cardinality of a set ${\sf S}^{\p}\left(
{\bf d}^3\right)$ \cite{frob87}, $t\left({\sf S}\left({\bf d}^3\right)\right)=
\#\left\{{\sf S}^{\p}\left({\bf d}^3\right)\right\}$, where 
\begin{eqnarray}
{\sf S}^{\p}\left({\bf d}^3\right)=\left\{x\in {\mathbb Z}\;\bracevert\;x
\not\in {\sf S}\left({\bf d}^3\right),x+s\in {\sf S}\left({\bf d}^3\right),\;
{\rm for\;all}\;s\in {\sf S}\left({\bf d}^3\right)/\{0\}\right\}.\nonumber
\end{eqnarray}
A set ${\sf S}^{\p}\left({\bf d}^3\right)$ is not empty since $F\left({\bf d}^
3\right)\in {\sf S}^{\p}\left({\bf d}^3\right)$ holds for any minimal 
generating set $\left(d_1,d_2,d_3\right)$.

Henceforth, ${\sf k}\left[{\sf S}\left({\bf d}^3\right)\right]$ is a 1--dim 
local Cohen--Macaulay ring of multiplicity $d_1$ and type $t\left({\sf S}\left(
{\bf d}^3\right)\right)$ which satisfies \cite{frob87}
\begin{eqnarray}
t\left({\sf S}\left({\bf d}^3\right)\right)=\left\{\begin{array}{l}1\;,\;\mbox{
if ${\sf S}\left({\bf d}^3\right)$ is symmetric}\;,\\2\;,\;\mbox{if ${\sf S}
\left({\bf d}^3\right)$ is non-symmetric}\;.\end{array}\right.\label{intro6}
\end{eqnarray}
Theorem \ref{theo1} and \ref{theo2} determine important relations 
between $G\left({\bf d}^3\right)$, $\widetilde{G}\left({\bf d}^3\right)$ and 
$t\left({\sf S}\left({\bf d}^3\right)\right)$.
\begin{theorem}{\rm (Theorem 20, \cite{frob87})}\label{theo1}
\begin{eqnarray}
G\left({\bf d}^3\right)\leq\widetilde{G}\left({\bf d}^3\right)
t\left({\sf S}\left({\bf d}^3\right)\right)\;.\label{intro7}
\end{eqnarray}
\end{theorem}
\begin{theorem}{\rm (Theorem 2, \cite{brown91}, Corollary at p. 339, 
\cite{brown92})}\label{theo2}
\begin{eqnarray}
G\left({\bf d}^3\right)=\left\{\begin{array}{l}\widetilde{G}\left({\bf d}^3
\right)\;,\;\;\mbox{iff}\;\;\;{\sf S}\left({\bf d}^3\right)\;\;\mbox{is
symmetric}\;,\\2\widetilde{G}\left({\bf d}^3\right)\;,\;\;\mbox{iff}\;\;\;{\bf
d}^3=\{3,3k+1,3k+2\}\;,\;\;k\geq 1\;.\end{array}\right.\label{intro8}
\end{eqnarray}
\end{theorem}
The rest of triples ${\bf d}^3$ gives rise to non--symmetric semigroups which 
satisfy a relation \cite{brown91},
\begin{eqnarray}
G\left({\bf d}^3\right)=2\widetilde{G}\left({\bf d}^3\right)-u\left({\bf d}^3
\right)\;,\;\;\;1\leq u\left({\bf d}^3\right)<\widetilde{G}\left({\bf d}^3 
\right)\;.\label{intro9}
\end{eqnarray}
The case $u\left({\bf d}^3\right)=1$ was studied in \cite{brown91}: it holds 
iff ${\sf S}\left({\bf d}^3\right)$ is generated by one of two sporadic triples
${\bf d}^3=\{4,5,11\},\;\{4,7,13\}$ or by one serie, ${\bf d}^3=\{3,3k+2,3k+4\},
\;k\geq 1$. As $u\left({\bf d}^3\right)$ increases, the number of sporadic 
triples climbs significantly. But there are not to our knowledge any general 
classification of such semigroups. However, it turns out that Theorems 
\ref{theo1} and \ref{theo2} are enough to resolve one of the Arnol'd Conjectures
(Conjecture \ref{con20}, see Section \ref{ccon2}). 

Consider a minimal generating set $\left(d_1,d_2,d_3\right)$ and let $g_1=
\gcd(d_2,d_3)$, $g_2=\gcd(d_3,d_1)$ and $g_3=\gcd(d_1,d_2)$ be given. We call 
the semigroup $\overline{{\sf S}}\left(d_1,d_2,d_3\right)$, 
\begin{eqnarray}
\overline{{\sf S}}\left(d_1,d_2,d_3\right)={\sf S}\left(\frac{d_1}{g_2g_3},
\frac{d_2}{g_1g_3},\frac{d_3}{g_1g_2}\right)\;,\label{wat1}
\end{eqnarray}
{\em the derived semigroup of ${\sf S}\left(d_1,d_2,d_3\right)$}.
\begin{theorem}{\rm (Corollary at p. 77, \cite{frob87})}\label{theo3}  
The semigroup ${\sf S}\left(d_1,d_2,d_3\right)$ is symmetric iff its derived 
semigroup $\overline{{\sf S}}\left(d_1,d_2,d_3\right)$ is generated by two 
elements.
\end{theorem}   
Now we state Theorem about necessary conditions for ${\sf S}\left({\bf d}^3
\right)$ to be symmetric.
\begin{theorem}\label{theo4}
If a semigroup ${\sf S}\left(d_1,d_2,d_3\right)$ is symmetric then its minimal
generating set has a following presentation with at least two relatively not 
prime elements:
\begin{eqnarray}
\gcd(d_1,d_2)=b\;,\;\;\gcd(d_3,b)=1\;,\;\;d_3\in {\sf S}\left(\frac{d_1}{b},
\frac{d_2}{b}\right)\;.\label{wat2}
\end{eqnarray}
\end{theorem}
{\sf Proof} $\;\;\;$ Let ${\sf S}\left(d_1,d_2,d_3\right)$ be a symmetric 
semigroup, i.e. $\gcd(d_1,d_2,d_3)=1$ and $\left(d_1,d_2,d_3\right)$ is a 
minimal generating set. According to Theorem \ref{theo3} its derived semigroup 
$\overline{{\sf S}}\left(d_1,d_2,d_3\right)$ given in (\ref{wat1}) is generated 
by two elements. Without loss of generality we can put $g_1=g_2=1$ and write,
\begin{eqnarray}
d_3=c_1\frac{d_1}{g_3}+c_2\frac{d_2}{g_3}\;,\;\;\;c_1,c_2\in {\mathbb Z}_+\;,
\label{wat3}
\end{eqnarray}
that results in
\begin{eqnarray}
\gcd(d_1,d_2)=g_3\;,\;\;\gcd(d_3,g_3)=1\;,\;\;d_3\in {\sf S}\left(\frac{d_1}
{g_3},\frac{d_2}{g_3}\right)\;.\label{wat4}
\end{eqnarray}
Denoting $g_3=b$ we arrive at (\ref{wat2}).$\;\;\;\;\;\;\Box$

It appears that (\ref{wat2}) gives also efficient conditions for ${\sf S}\left(
{\bf d}^3\right)$ to be symmetric. This follows from Corollary of the early 
Lemma \cite{wata73} for semigroup ${\sf S}\left({\bf d}^m\right)$
\begin{lemma}{\rm (Lemma 1, \cite{wata73})}\label{lem1}   
Let ${\sf S}\left(d_1,\ldots,d_m\right)$ be a numerical semigroup, $a$ and $b$ 
be positive integers such that: (i) $a\in {\sf S}\left(d_1,\ldots,d_m\right)$
and $a\neq d_i$, (ii) $\gcd(a,b)=1$.\\
Then a semigroup ${\sf S}\left(bd_1,\ldots,bd_m,a\right)$ is symmetric iff 
${\sf S}\left(d_1,\ldots,d_m\right)$ is symmetric.
\end{lemma}
Combining Lemma \ref{lem1} with a fact that a semigroup ${\sf S}\left({\bf 
d}^2\right)$ is always symmetric we arrive at Corollary.
\begin{corollary}\label{cor1}  
Let ${\sf S}\left(d_1,d_2\right)$ be a numerical semigroup, $a$ and $b$ be
positive integers, $\gcd(a,b)=1$. If $a\in {\sf S}\left(d_1,d_2\right)$, then 
a semigroup ${\sf S}\left(bd_1,bd_2,a\right)$ is symmetric.
\end{corollary}
In Corollary \ref{cor1} a requirement $a\neq d_1,d_2$ can be omitted since 
e.g. a semigroup ${\sf S}\left(bd_1,bd_2,d_1\right)$ is generated by two 
elements ($d_1,bd_2$) and is also symmetric. 

For a sake of completeness finish this Section with efficient and necessary 
conditions for ${\sf S}\left({\bf d}^3\right)$ to be non--symmetric.
\begin{theorem}{\rm (Theorem 14, \cite{frob87})}\label{theo5}
A semigroup ${\sf S}\left(d_1,d_2,d_3\right)$ minimally generated by three 
pairwise relatively prime elements is non--symmetric.
\end{theorem}
\begin{theorem}{\rm (Corollary at p. 71, \cite{frob87})}\label{theo6}
Let ${\sf S}\left(d_1,d_2,d_3\right)$ be a semigroup and $\overline{{\sf S}}
\left(d_1,d_2,d_3\right)$ be its derived semigroup. Then ${\sf S}\left(d_1,d_2,
d_3\right)$ and $\overline{{\sf S}}\left(d_1,d_2,d_3\right)$ have the same
type. In particular, ${\sf S}\left(d_1,d_2,d_3\right)$ is non--symmetric iff
$\overline{{\sf S}}\left(d_1,d_2,d_3\right)$ is non--symmetric. 
\end{theorem}
More specific details on semigroups ${\sf S}\left({\bf d}^3\right)$ will be 
given in Section \ref{stat0}.
\section{Weak asymptotics in numerical semigroups ${\sf S}\left({\bf d}^m
\right)$}\label{stat1}
Two sequences of real numbers $A(k)$ and $B(k),\;k\in {\mathbb Z}_+$, are
said to have the same {\em weak asymptotics} \cite{arn99}, or to have the
same {\em growth rate} \cite{arn04}, or to be {\em Ces\'aro equivalent}
\cite{arn06}, if
\begin{eqnarray}
\lim_{N\to \infty}\frac{\sum_{k=1}^NA(k)}{\sum_{k=1}^NB(k)}=1\;.\label{res0}
\end{eqnarray}
The limit is weak here: one requires the convergence only for the sums in
(\ref{res0}). In the similar way one can consider the Ces\'aro equivalence of 
two sequences $A(k)$ and $B(k)$ {\em at large integers} $k\in {\mathbb Z}_+$. 
Let us replace $k$ by a neighborhood ${\mathbb U}_{N,r}(k)$ of length $2r$ of 
a scaled integer $Nk,N\in {\mathbb Z}_+$. Replace the values of $A(k)$ and 
$B(k)$ by the arithmetic means $A_{N,r}(k)$ and $B_{N,r}(k)$, respectively,
\begin{eqnarray}
A_{N,r}(k)=\frac{1}{2r}\sum_{j=-r}^{r}A(Nk+j)\;,\;\;\;\;B_{N,r}(k)=\frac{1}
{2r}\sum_{j=-r}^{r}B(Nk+j)\;,\;\;\;\;Nk+j\in {\mathbb U}_{N,r}(k)\;.\nonumber
\end{eqnarray}
Two sequences of real numbers $A(k)$ and $B(k),\;k\in {\mathbb Z}_+$, are 
said to have the same {\em weak asymptotics at large} $k$ \cite{arn06}, if
\begin{eqnarray}
\lim_{r,N\to\infty\atop r(N)/N\to 0}\frac{A_{N,r}(k)}{B_{N,r}(k)}=1\;.
\label{res1}
\end{eqnarray}
A study of weak asymptotics {\em at the typical large vectors} ${\bf d}^m$ for 
numerical functions $A\left({\bf d}^m\right)$, as conductor or genus of 
semigroup, over all vectors ${\bf b}^m$ comprising a numerical semigroup ${\sf 
S}\left({\bf d}^m\right)$ is much more difficult problem. Call such vectors 
${\bf b}^m$, ${\bf b}^m\in {\sf S}\left({\bf d}^m\right)$, {\em typical}. 
Arnol'd gave a recept \cite{arn04} how to average such functions over {\em the 
typical large vectors} ${\bf b}^m$.

Let ${\sf S}\left({\bf d}^m\right)$ be a numerical semigroup, i.e. a generating 
set $\left(d_1,\ldots,d_m\right)$ is minimal. Replace the vector ${\bf d}^m$ 
by a spheric (or cubic) neighborhood ${\mathbb U}_{N,r}\left({\bf d}^m\right)$ 
of radius $r$ of a scaled vector $N{\bf d}^m\in {\mathbb Z}_+^m,\;N\in
{\mathbb Z}_+$. Denote by ${\bf j}^m$ a 
vector $(j_1,\ldots,j_m)$. Replace the value $A\left({\bf d}^m\right)$ by the 
arithmetic mean $A_{N,r}\left({\bf d}^m\right)$ of the functions $A\left(N{\bf 
d}^m+{\bf j}^m\right)$ at the vectors $N{\bf d}^m+{\bf j}^m\in {\mathbb U}_{
N,r}\left({\bf d}^m\right)$ whose components $Nd_i+j_i,j_i\in {\mathbb Z}_+,\;
-r\leq j_i\leq r$, satisfy two constraints:
\begin{enumerate}
\item The following holds,
\vspace{-.4cm}
\begin{equation}
\gcd(Nd_1+j_1,\ldots,Nd_m+j_m)=1\;,\label{resc2}
\end{equation}
otherwise the corresponding numerical semigroup has an infinite complement 
$\Delta\left({\bf d}^m\right)$. 
\item $\{Nd_1+j_1,\ldots,Nd_m+j_m\}$ is {\em a minimal generating set}, i.e.
there are no nonnegative integers $f_{i,k}$ for which a linear dependence holds
\begin{equation} 
Nd_i+j_i=\sum_{k\neq i}^mf_{i,k}(Nd_k+j_k)\;,\;\;\;f_{i,k}\in\{0,1,\ldots\}
\;\;\;\mbox{for any}\;\;i\leq m\;,\label{resc1}
\end{equation}  
otherwise $N{\bf d}^m+{\bf j}^m$ does not generate the m--dim numerical
semigroup.
\end{enumerate}
We have also choose the averaging radius $r$ and its growth rate such that
\begin{eqnarray}
1\ll r\ll N\;,\;\;\;r(N)/N\to 0\;\;\;\mbox{when}\;\;\;N\to\infty\;.\label{s00c}
\end{eqnarray}
Call the vector $N{\bf d}^m+{\bf j}^m$ {\em admissible} if its components
satisfy both constraints (\ref{resc2}) and (\ref{resc1}).  Denote by ${\mathbb 
M}_{N,r}\left({\bf d}^m\right)$ an entire set of admissible vectors, ${\mathbb 
M}_{N,r}\left({\bf d}^m\right)\subset {\mathbb U}_{N,r}\left({\bf d}^m\right)$, 
\begin{equation}
{\mathbb M}_{N,r}\left({\bf d}^m\right)=\left\{N{\bf d}^m+{\bf j}^m|-r\leq j_i
\leq r,1\ll r\ll N,{\rm Constraints\;(\ref{resc2})\;and\;(\ref{resc1})\;are\;
satisfied}\right\}.\label{set1}
\end{equation}
Denote by $\#\left\{{\mathbb M}_{N,r}\left({\bf d}^m\right)\right\}$ a 
cardinality of ${\mathbb M}_{N,r}\left({\bf d}^m\right)$ and notice that $\#
\left\{{\mathbb M}_{N,r}\left({\bf d}^m\right)\right\}<(2r)^m$ since at least 
$N{\bf d}^m\not\in {\mathbb M}_{N,r}\left({\bf d}^m\right)$ because $\gcd(Nd_1,
\ldots,Nd_m)=N$. Write the arithmetic mean, 
\begin{eqnarray}  
A_{N,r}\left({\bf d}^m\right)=\frac{1}{\#\left\{{\mathbb M}_{N,r}\left({\bf 
d}^m\right)\right\}}\sum_{j_1,\ldots,j_m=-r}^{r}A\left(N{\bf d}^m+{\bf j}^m
\right)\;,\;\;\;\;N{\bf d}^m+{\bf j}^m\in {\mathbb M}_{N,r}\left({\bf d}^m
\right)\;.\label{set1a}
\end{eqnarray}
Say that two numerical functions $A\left({\bf d}^m\right)$ and $B\left({\bf d}
^m\right)$ have the same {\em weak asymptotics at the typical large} ${\bf d}^
m$ \cite{arn04}, if   
\begin{eqnarray}
\lim_{r,N\to\infty\atop r(N)/N\to 0}\frac{A_{N,r}\left({\bf d}^m\right)}{
B_{N,r}\left({\bf d}^m\right)}=\lim_{r,N\to\infty\atop r(N)/N\to 0}\frac{
\sum_{j_1,\ldots,j_m=-r}^{r}A\left(N{\bf d}^m+{\bf j}^m\right)}{\sum_{j_1,
\ldots,j_m=-r}^{r}B\left(N{\bf d}^m+{\bf j}^m\right)}=1\;,\;\;\;N{\bf d}^m+
{\bf j}^m\in {\mathbb M}_{N,r}\left({\bf d}^m\right)\;,\label{set1b}
\end{eqnarray}
and denote this equivalence by
\begin{eqnarray}
{\sf A}\left({\bf d}^m\right)\stackrel{asymptotically\atop weak}{\equiv}
{\sf B}\left({\bf d}^m\right)\;.\label{set1c}
\end{eqnarray}
\section{Statistics of numerical semigroups ${\sf S}\left(N{\bf d}^m+{\bf j}^m
\right),\;N\to\infty$}\label{numer1}
The main difficulties in performing an analytic summation in (\ref{set1a}) and 
(\ref{set1b}) are caused by constraints (\ref{resc2}) and (\ref{resc1}) which 
are hardly to account for. For this aim let us estimate $\#\left\{{\mathbb M}_
{N,r}\left({\bf d}^m\right)\right\}$ in the limit (\ref{s00c}). 
Represent a set ${\mathbb M}_{N,r}\left({\bf d}^m\right)$ as follows, 
\begin{eqnarray}  
{\mathbb M}_{N,r}\left({\bf d}^m\right)&=&\widehat{{\mathbb M}_{N,r}}\left({\bf 
d}^m\right)\setminus \widetilde{{\mathbb M}_{N,r}}\left({\bf d}^m\right)\;,
\label{set1c}\\
\widehat{{\mathbb M}_{N,r}}\left({\bf d}^m\right)&=&\left\{N{\bf d}^m+{\bf 
j}^m\;|\;-r\leq j_i\leq r,\;1\ll r\ll N,\;{\rm Constraint\;(\ref{resc2})\;is
\;satisfied}\right\}\;,\label{set1d}
\end{eqnarray}
and a set $\widetilde{{\mathbb M}_{N,r}}\left({\bf d}^m\right)$ comprises all 
vectors $N{\bf d}^m+{\bf j}^m$ whose generating sets $\{Nd_1+j_1,\ldots,Nd_m+j_
m\}$ are not minimal though they are still satisfying Constraint (\ref{resc2}). 
Consider two sets $\widehat{{\mathbb M}_{N,r}}\left({\bf d}^m\right)$ and
$\widetilde{{\mathbb M}_{N,r}}\left({\bf d}^m\right)$ separately.

Calculate a cardinality of a set $\widehat{{\mathbb M}_{N,r}}\left({\bf d}^m
\right)$ in the limit (\ref{s00c}) by a probabilistic method which dates back 
to Euler \cite{arn06}. It is based on a geometric interpretation of probability 
${\cal P}_{r,\infty}$ that randomly chosen integers $(Nd_1+j_1,\ldots,Nd_m+
j_m)$ from a set ${\mathbb U}_{N,r}\left({\bf d}^m\right)$ do not have common 
divisors
\begin{eqnarray}
{\cal P}_{r,\infty}=\lim_{N\to \infty}\frac{\#\left\{\widehat{{\mathbb M}_{N,r}
}\left({\bf d}^m\right)\right\}}{\#\left\{{\mathbb U}_{N,r}\left({\bf d}^m
\right)\right\}}=\frac{1}{(2r)^m}\lim_{N\to \infty}\#\left\{\widehat{{\mathbb 
M}_{N,r}}\left({\bf d}^m\right)\right\}\;.\label{set1j}
\end{eqnarray}
Let a tuple $\{Nd_1+j_1,\ldots,Nd_m+j_m\}$ is chosen
randomly from a cubic neighborhood ${\mathbb U}_{N,r}\left({\bf d}^m\right)$ 
of a scaled vector $N{\bf d}^m\in {\mathbb Z}_+^m,\;N\in{\mathbb Z}_+$ with the 
edge length $2r$ such that $1\ll r\ll N$. The least integer which is still 
contained in ${\mathbb U}_{N,r}\left({\bf d}^m\right)$ is $Nd_1-r$. Let $p$ be 
a prime integer such that $p\leq Nd_1-r$. A probability that $p$ divides every 
element $Nd_i+j_i$ in a tuple $\{Nd_1+j_1,\ldots,Nd_m+j_m\}$ is given by $p^{-
m}$. Consequently, $1-p^{-m}$ is a probability that $p$ does not divide any 
element in this tuple. Multiplying it over all primes such that $p\leq Nd_1-r$, 
we arrive
\begin{eqnarray}
{\cal P}_{r,N}=\prod_{2\leq p\leq Nd_1-r}\left(1-\frac{1}{p^m}\right)\;,
\label{set1f}
\end{eqnarray}
where ${\cal P}_{r,N}$ gives a probability that integers $(Nd_1+j_1,\ldots,
Nd_m+j_m)$, which are randomly chosen from a set ${\mathbb U}_{N,r}\left({\bf 
d}^m\right)$, do not have common divisors in the range $2\leq p\leq Nd_1-r$. 
Taking the limit (\ref{s00c}) we get
\begin{eqnarray}
{\cal P}_{r,\infty}=\prod_{p}^{\infty}\left(1-\frac{1}{p^m}\right)=
\frac{1}{\zeta(m)}\;,\label{set1g}
\end{eqnarray}
where $\zeta(m)$ stands for the Riemann zeta function. The value $\zeta^{-1}
(m)$ gives a probability that there are no other integral points on the segment 
between 0 and an integral point in $m$--dim space \cite{arn06}. Its first 
several values read $\zeta^{-1}(2)=0.6079$, $\zeta^{-1}(3)=0.8319$, $\zeta^{-1}
(4)=0.9239$. 

Thus, a cardinality of a set $\widehat{{\mathbb M}_{N,r}}\left({
\bf d}^m\right)$ can be estimated as 
\begin{eqnarray}
\#\left\{\widehat{{\mathbb M}_{N,r}}\left({\bf d}^m\right)\right\}\simeq 
\frac{(2r)^m}{\zeta(m)}\;.\label{set1e}
\end{eqnarray}
As for a set $\widetilde{{\mathbb M}_{N,r}}\left({\bf d}^m\right)$, the 
constraint (\ref{resc2}) in the case $m=2$ already presumes that $\widehat{{
\mathbb M}_{N,r}}\left({\bf d}^2\right)$ comprises all vectors $N{\bf d}^2+
{\bf j}^2$ whose generating sets $\{Nd_1+j_1,Nd_2+j_2\}$ are minimal, and 
therefore
\begin{eqnarray}
\widetilde{{\mathbb M}_{N,r}}\left({\bf d}^2\right)=\emptyset\;,\;\;\;
{\mathbb M}_{N,r}\left({\bf d}^2\right)&=&\widehat{{\mathbb M}_{N,r}}\left(
{\bf d}^2\right)\;,\;\;\;\#\left\{{\mathbb M}_{N,r}\left({\bf d}^2\right)
\right\}\simeq\frac{4r^2}{\zeta(2)}\;.\label{set1i} 
\end{eqnarray}
Calculate a cardinality of a set $\widetilde{{\mathbb M}_{N,r}}\left({\bf d}^m
\right)$ in the limit (\ref{s00c}) for higher $m$, $m\geq 3$. Let a vector 
$N{\bf d}^m+{\bf j}^m\in \widetilde{{\mathbb M}_{N,r}}\left({\bf d}^m\right)$  
be given, i.e. a generating set $\{Nd_1+j_1,\ldots,Nd_m+j_m\}$ is not minimal. 
According to (\ref{resc1}) there exists at least one element $Nd_i+j_i$ which 
is representable through the rest of the tuple,
\begin{eqnarray}
Nd_i+j_i=\sum_{k\neq i}^mf_{i,k}(Nd_k+j_k)\;,\;\;\;\;\;\mbox{or}\;\;\;\;\;d_i-
\sum_{k\neq i}^mf_{i,k}d_k=\frac{1}{N}\left(\sum_{k\neq i}^mf_{i,k}j_k-j_i
\right)\;,\label{set1k}
\end{eqnarray}
where $f_{i,k}\in\{0,1,\ldots\}$. Taking the limit (\ref{s00c}) we get two 
relations imposed on $N{\bf d}^m+{\bf j}^m$, 
\begin{eqnarray}
d_i=\sum_{k\neq i}^mf_{i,k}d_k\;,\;\;\;\;\;\mbox{and}\;\;\;\;\;
j_i=\sum_{k\neq i}^mf_{i,k}j_k\;.\label{set1l}
\end{eqnarray}
The 1st of them claims that the generating set $\{d_1,\ldots,d_m\}$ is not 
minimal. However, this contradicts an assumption that ${\sf S}\left({\bf d}^m
\right)$ is a numerical semigroup generated by $m$ elements. Thus, a relation 
(\ref{set1k}) can not be satisfied by any choice of $(j_1,\ldots,j_m\}$ and 
therefore a set $\widetilde{{\mathbb M}_{N,r}}\left({\bf d}^m\right)$ in the 
limit (\ref{s00c}) doesn't contain any vectors. Thus, we have
\begin{eqnarray}
\widetilde{{\mathbb M}_{N,r}}\left({\bf d}^m\right)=\emptyset\;,\;\;\;
{\mathbb M}_{N,r}\left({\bf d}^m\right)&=&\widehat{{\mathbb M}_{N,r}}\left(
{\bf d}^m\right)\;,\;\;\;\#\left\{{\mathbb M}_{N,r}\left({\bf d}^m\right)
\right\}\simeq\frac{(2r)^m}{\zeta(m)}\;.\label{set1m}
\end{eqnarray}
\subsection{Statistics of symmetric and non-symmetric semigroups ${\sf S}
\left(N{\bf d}^3+{\bf j}^3\right),\;N\to\infty$}\label{met1}
In this Section we deal with numerical semigroups ${\sf S}\left(N{\bf d}^3+
{\bf j}^3\right)$ generated by three elements only. 
Consider statistics of symmetric and non--symmetric semigroups ${\sf S}\left(
N{\bf d}^3+{\bf j}^3\right)$ corresponding to admissible vectors $N{\bf d}^3+
{\bf j}^3$. Denote by ${\mathbb M}_{N,r}^{sym}\left({\bf d}^3\right)$ and 
${\mathbb M}_{N,r}^{nsym}\left({\bf d}^3\right)$ the sets of admissible vectors 
$N{\bf d}^3+{\bf j}^3\in{\mathbb M}_{N,r}\left({\bf d}^3\right)$ such that 
they correspond to the minimal generating sets of symmetric and non--symmetric 
semigroups, respectively,
\begin{eqnarray}
{\mathbb M}_{N,r}^{sym}\left({\bf d}^3\right)&=&\left\{N{\bf d}^3+{\bf j}^3|-r
\leq j_i\leq r,1\ll r\ll N,\;{\sf S}\left(N{\bf d}^3+{\bf j}^3\right)\;{\rm\;
is\;symmetric}\right\},\label{kag1}\\
{\mathbb M}_{N,r}^{nsym}\left({\bf d}^3\right)&=&\left\{N{\bf d}^3+{\bf j}^3|-r
\leq j_i\leq r,1\ll r\ll N,\;{\sf S}\left(N{\bf d}^3+{\bf j}^3\right)\;{\rm\;
is\;non-symmetric}\right\}.\nonumber
\end{eqnarray}
These sets and their cardinalities ($\#$) are related in the following way,
\begin{eqnarray}
{\mathbb M}_{N,r}\left({\bf d}^3\right)&=&{\mathbb M}_{N,r}^{sym}\left({\bf d}^
3\right)\;\cup\;{\mathbb M}_{N,r}^{nsym}\left({\bf d}^3\right)\;,\nonumber\\
\#\left\{{\mathbb M}_{N,r}\left({\bf d}^3\right)\right\}&=&\#\left\{{\mathbb M}
_{N,r}^{sym}\left({\bf d}^3\right)\right\}+\#\left\{{\mathbb M}_{N,r}^{nsym}
\left({\bf d}^3\right)\right\}\;.\label{kag2}
\end{eqnarray}
Calculate a cardinality of a set ${\mathbb M}_{N,r}^{sym}\left({\bf d}^3
\right)$ in the limit (\ref{s00c}) by applying Theorem \ref{theo4}. 
\begin{lemma}\label{lem2}
Let ${\sf S}\left({\bf d}^3\right)$ and ${\sf S}\left(N{\bf d}^3+{\bf j}^3
\right)$ be numerical semigroups and $\{Nd_1+j_1,Nd_2+j_2,Nd_3+j_3\}$ be 
a minimal generating set such that
\begin{eqnarray}
-r\leq j_1,j_2,j_3\leq r\;,\;\;1\ll r\ll N\;.\label{kag2a}
\end{eqnarray}
If $r(N)/N\to 0$ when $N\to\infty$ then
\begin{eqnarray}
\lim_{N\to 0}\#\left\{{\mathbb M}_{N,r}^{sym}\left({\bf d}^3\right)\right\}=0
\;.\label{kag2b}
\end{eqnarray}
\end{lemma}
{\sf Proof} $\;\;\;$ Consider a minimal generating set $\{Nd_1+j_1,Nd_2+j_2,
Nd_3+j_3\}$ satisfying (\ref{kag2a}). Suppose that the corresponding semigroup 
${\sf S}\left(N{\bf d}^3+{\bf j}^3\right)$ is symmetric. According to Theorem 
\ref{theo4} a triple $N{\bf d}^3+{\bf j}^3$ has necessarily a following 
presentation:
\begin{eqnarray}
&&\gcd\left(Nd_1+j_1,Nd_2+j_2\right)=b\;,\;\;\;b\in {\mathbb Z}_+\;,\;\;
b\geq 2\;,\label{kag3}\\
&&\gcd\left(Nd_3+j_3,b\right)=1\;,\label{kag3a}\\
&&b\left(Nd_3+j_3\right)=c_1\left(Nd_1+j_1\right)+c_2\left(Nd_2+j_2\right)\;,
\;\;\;c_1,c_2\in {\mathbb Z}_+\;.\label{kag4}
\end{eqnarray}
First, consider (\ref{kag3}) and find a great common divisor $b$ of integers 
$Nd_1+j_1$ and $Nd_2+j_2$ in the limit (\ref{s00c}). Let such $b$ exists, then
\begin{eqnarray}
Nd_1+j_1=b\;k_1\;,\;\;\;Nd_2+j_2=b\;k_2\;,\;\;\;k_1,k_2\in {\mathbb Z}_+\;,
\;\;\;\gcd\left(k_1,k_2\right)=1\;,\label{kag4a}
\end{eqnarray}
that results in the following
\begin{eqnarray}
k_1(Nd_2+j_2)=k_2(Nd_1+j_1)\;,\;\;\;\rightarrow\;\;\;k_1d_2-k_2d_1=\frac{1}
{N}\left(k_2j_1-k_1j_2\right)\;.\label{kag4b}
\end{eqnarray}
Taking the limit (\ref{s00c}) we get two Diophantine equations 
\begin{eqnarray}
k_1d_2-k_2d_1=0\;,\;\;\;\;\;k_1j_2-k_2j_1=0\;,\label{kag4y}
\end{eqnarray}
supplemented by $\gcd\left(k_1,k_2\right)=1$. Their solutions read
\begin{eqnarray}
k_1=\frac{{\sf lcm}(d_1,d_2)}{d_2}\;,\;\;k_2=\frac{{\sf lcm}(d_1,d_2)}{d_1}\;,
\;\;\;\;j_1=k_3d_1\;,\;\;j_2=k_3d_2\;,\;\;k_3\in {\mathbb Z}_+\;,\;\;
k_3\leq \frac{r}{d_2}\;.\label{kag4c}
\end{eqnarray}
Combining (\ref{kag4a}) and (\ref{kag4c}) we get
\begin{eqnarray}
b=(N+k_3)\gcd(d_1,d_2)\;.\label{kag4d}
\end{eqnarray}
As for (\ref{kag3a}), its simple comparison with (\ref{kag4d}) necessarily 
claims
\begin{eqnarray}
j_3\neq k_3d_3\;.\label{kag4e}
\end{eqnarray}
Finally, consider a relation (\ref{kag4}) in the limit (\ref{s00c}). Similarly 
to (\ref{kag4y}) we get two Diophantine equations imposed on the tuples 
$\left(d_1,d_2,d_3\right)$ and $\left(j_1,j_2,j_3\right)$ separately,
\begin{eqnarray}
bd_3=c_1d_1+c_2d_2\;,\;\;\;\;\;\mbox{and}\;\;\;\;\;bj_3=c_1j_1+c_2j_2\;.
\label{kag5}
\end{eqnarray}
Multiplying the 1st equation by $k_3$ and making difference between both  
equations we get
\begin{eqnarray}
b(j_3-k_3d_3)=0\;,\label{kag5a}
\end{eqnarray}
that contradicts (\ref{kag4e}). Thus, a set ${\mathbb M}_{N,r}^{sym}\left({\bf 
d}^3\right)$ is empty in the limit (\ref{s00c}) that proves Lemma.
$\;\;\;\;\;\;\Box$

Thus, by (\ref{set1m}) and (\ref{kag2}) we have
\begin{eqnarray}
\#\left\{{\mathbb M}_{N,r}^{nsym}\left({\bf d}^3\right)\right\}\simeq
\#\left\{{\mathbb M}_{N,r}\left({\bf d}^3\right)\right\}\simeq\frac{(2r)^3}
{\zeta(3)}\;.\label{kag6}
\end{eqnarray}
In other words, a summation in (\ref{set1a}) and (\ref{set1b}) for $m=3$ is 
performed over all non--symmetric semigroups ${\sf S}\left(N{\bf d}^3+{\bf j}^3
\right)$ exclusively. 
\section{Arnold's conjectures on weak asymptotics}\label{arr1}
V. Arnold gave his conjectures on weak asymptotics for the numerical semigroups
${\sf S}\left({\bf d}^m\right)$ of arbitrary dimension $m$.
\subsection{Conjecture $\#$1999--8 and its discussion}\label{ccon1}
Conjecture $\#$1999--8 deals with the asymptotic behavior of the conductor
of the numerical semigroups. We quote from \cite{arn00} :
\begin{conjecture}{\rm ($\#$1999--8)}\label{con10}
Explore the statistics of $C\left({\bf d}^m\right)$ for typical large vectors
${\bf d}^m$. \\Conjecturally,
\begin{eqnarray}
{\sf C}\left({\bf d}^m\right)\stackrel{asymptotically\atop weak}{\equiv}g_m
\sqrt[m-1]{d_1\cdot\ldots\cdot d_m}\;,\;\;\;\;\;\;g_m=\sqrt[m-1]{(m-1)!}\;.
\label{syl1}
\end{eqnarray}
\end{conjecture}
Define a new function $K_{N,r}\left({\bf d}^m\right)$ in the sense of 
(\ref{set1b}),
\begin{eqnarray}
K_{N,r}\left({\bf d}^m\right)=\frac{\sum_{j_1,\ldots,j_m=-r}^{r}C\left(N{\bf 
d}^m+{\bf j}^m\right)}{\sum_{j_1,\ldots,j_m=-r}^{r}\sqrt[m-1]{V\left(N{\bf d}^
m+{\bf j}^m\right)}}\;,\;\;\;\mbox{where}\;\;\;\;N{\bf d}^m+{\bf j}^ m\in 
{\mathbb M}_{N,r}\left({\bf d}^m\right)\;,\label{s0c}
\end{eqnarray}
and $V\left({\bf d}^m\right)=d_1\cdot\ldots\cdot d_m$. 
Then Conjecture \ref{con10} can been represented as follows,
\begin{eqnarray}
{\sf K}\left({\bf d}^m\right)=g_m\;,\;\;\;\mbox{where}\;\;\;{\sf K}\left({\bf 
d}^m\right)=\lim_{r,N\to\infty\atop r(N)/N\to 0}K_{N,r}\left({\bf d}^m\right)\;.
\label{s1c}
\end{eqnarray}
For $m=2$ the corresponding semigroups ${\sf S}\left(N{\bf d}^2+{\bf j}^2
\right)$ are symmetric and the problem is simplified essentially due to the two 
reasons. First, the constraint (\ref{resc1}) is already incorporated into 
(\ref{resc2}). Next, the conductor $C\left({\bf d}^2\right)$ is known due to 
Sylvester \cite{sylv84}, $C\left(d_1,d_2\right)=(d_1-1)(d_2-1)$. Performing the 
calculation we can verify Conjecture \ref{con10} for $m=2$,
\begin{eqnarray}
{\sf K}\left({\bf d}^2\right)=\lim_{r,N\to\infty\atop r(N)/N\to 0}\frac{\sum_
{j_1,j_2=-r\atop \gcd(N{\bf d}^2+{\bf j}^2)=1}^{r}\left(1+\frac{j_1-1}{Nd_1}
\right)\left(1+\frac{j_2-1}{Nd_2}\right)}{\sum_{j_1,j_2=-r\atop\gcd(N{\bf d}^2
+{\bf j}^2)=1}^{r}\left(1+\frac{j_1}{Nd_1}\right)\left(1+\frac{j_2}{Nd_2}
\right)}=1-\lim_{r,N\to\infty\atop r(N)/N\to 0}A_1+\lim_{r,N\to\infty\atop 
r(N)/N\to 0}A_2\;,\label{s2c}
\end{eqnarray}
where
\begin{eqnarray}
A_1&\simeq&\frac{r}{N}\;\frac{d_1^{-1}\sum_{j_2=-r\atop \gcd(N{\bf d}^2+{\bf j}
^2)=1}^{r}\left(1+\frac{j_2}{Nd_2}\right)+d_2^{-1}\sum_{j_1=-r\atop\gcd(N{\bf 
d}^2+{\bf j}^2)=1}^{r}\left(1+\frac{j_1}{Nd_2}\right)}{\sum_{j_1,j_2=-r\atop
\gcd(N{\bf d}^2+{\bf j}^2)=1}^{r}\left(1+\frac{j_1}{Nd_1}\right)\left(1+\frac{
j_2}{Nd_2}\right)}\;,\nonumber\\
A_2&\simeq&\left(\frac{r}{N}\right)^2\frac{d_1^{-1}d_2^{-1}}{\sum_{j_1,j_2=-r
\atop\gcd(N{\bf d}^2+{\bf j}^2)=1}^{r}\left(1+\frac{j_1}{Nd_1}\right)\left(1+
\frac{j_2}{Nd_2}\right)}\;.\nonumber
\end{eqnarray}
Taking the limit (\ref{s00c}) in (\ref{s2c}) we have ${\sf K}\left({\bf d}^2
\right)=1$. 

For $m\geq 3$ the main difficulty in performing an analytic summation in 
(\ref{s1c}) is due to Curtis' theorem \cite{curt90} on the non--algebraic 
representation of the Frobenius number $F\left({\bf d}^m\right)$. In other 
words, $F\left({\bf d}^m\right)$ cannot be expressed by $d_1,\ldots,d_m$ as 
an algebraic function (see also \cite{fel05}). In order to overcome this 
difficulty and discuss Conjecture \ref{con10} in the case $m=3$ we will bound 
the limit in (\ref{s1c}).

Consider the 3--dim version of Conjecture \ref{con10} and recall recent 
results \cite{fel04} about the lower bounds for conductor in the symmetric 
and non--symmetric semigroups ${\sf S}\left({\bf d}^3\right)$,
\begin{eqnarray}
C\left({\bf d}^3\right)\geq\left\{\begin{array}{l}\sqrt{3}\sqrt{d_1d_2d_3+1}-
(d_1+d_2+d_3)+1\;,\;\mbox{if ${\sf S}\left({\bf d}^3\right)$ is non-symmetric}
\;,\\
2\sqrt{d_1d_2d_3}-(d_1+d_2+d_3)+1\;,\;\mbox{if ${\sf S}\left({\bf d}^3\right)$ 
is symmetric}\;.\end{array}\right.\label{lllow1}
\end{eqnarray}
Define a ratio,
\begin{eqnarray}
v\left(N{\bf d}^3+{\bf j}^3\right)=\frac{C\left(N{\bf d}^3+{\bf j}^3\right)}
{\sqrt{V\left(N{\bf d}^3+{\bf j}^3\right)}}\;,\label{s3c1}
\end{eqnarray}
and represent $K_{N,r}\left({\bf d}^3\right)$ as follows,
\begin{eqnarray}
K_{N,r}\left({\bf d}^3\right)=\frac{\sum_{j_1,j_2,j_3=-r\atop
\gcd(N{\bf d}^3+{\bf j}^3)=1}^{r}v\left(N{\bf d}^3+{\bf j}^3\right)\sqrt{V
\left(N{\bf d}^3+{\bf j}^3\right)}}{\sum_{j_1,j_2,j_3=-r\atop \gcd(N{\bf d}^3
+{\bf j}^3)=1}^{r}\sqrt{V\left(N{\bf d}^3+{\bf j}^3\right)}}\;\;.\label{s3c}
\end{eqnarray}
By Lemma \ref{lem2} a summation in (\ref{s3c}) is performed over $N{\bf d}^3+
{\bf j}^3\in {\mathbb M}_{N,r}\left({\bf d}^3\right)$ and the corresponding 
semigroups ${\sf S}\left(N{\bf d}^3+{\bf j}^3\right)$ are non--symmetric only.  
A bound (\ref{lllow1}) is valid for all such admissible vectors. This results 
in the following,
\begin{eqnarray}
v\left(N{\bf d}^3+{\bf j}^3\right)&\geq&\frac{\sqrt{3}\sqrt{(Nd_1+j_1)
(Nd_2+j_2)(Nd_3+j_3)+1}-N\sum_{i=1}^3d_i-\sum^3_{i=1}j_i+1}{\sqrt{(Nd_1+j_1)
(Nd_2+j_2)(Nd_3+j_3)}}>\nonumber\\
&&\sqrt{3}-\frac{d_1+d_2+d_3}{\sqrt{N}\sqrt{d_1d_2d_3}}\frac{1+
\frac{j_1+j_2+j_3-1}{N(d_1+d_2+d_3)}}{\sqrt{\left(1+\frac{j_1}{Nd_1}\right)
\left(1+\frac{j_2}{Nd_2}\right)\left(1+\frac{j_3}{Nd_3}\right)}}\;.\label{s3c2}
\end{eqnarray}
Combining (\ref{s3c}) and (\ref{s3c2}) we get
\begin{eqnarray}
K_{N,r}\left({\bf d}^3\right)>\sqrt{3}-\frac{d_1+d_2+d_3}{\sqrt{N}\sqrt{d_1d_2
d_3}}\left(\frac{\Sigma_1(r,N)}{\Sigma_0(r,N)}+\frac{\Sigma_2(r,N)}{\Sigma_0(r,
N)}\right)\;,\label{s3c3}
\end{eqnarray}
where
\begin{eqnarray}
\Sigma_0(r,N)&=&\sum_{j_1,j_2,j_3=-r\atop\gcd(N{\bf d}^3+{\bf j}^3)=1}^{r}
\sqrt{V\left(N{\bf d}^3+{\bf j}^3\right)}\nonumber\\
\Sigma_1(r,N)&=&\sum_{j_1,j_2,j_3=-r\atop\gcd(N{\bf d}^3+{\bf j}^3)=1}^{r}
\sqrt{\frac{V\left(N{\bf d}^3+{\bf j}^3\right)}{\left(1+\frac{j_1}{Nd_1}
\right)\left(1+\frac{j_2}{Nd_2}\right)\left(1+\frac{j_3}{Nd_3}\right)}}\;,
\nonumber\\
\Sigma_2(r,N)&=&\sum_{j_1,j_2,j_3=-r\atop\gcd(N{\bf d}^3+{\bf j}^3)=1}^{r}
\frac{j_1+j_2+j_3-1}{N(d_1+d_2+d_3)}\sqrt{\frac{V\left(N{\bf d}^3+{\bf j}^3
\right)}{\left(1+\frac{j_1}{Nd_1}\right)\left(1+\frac{j_2}{Nd_2}\right)
\left(1+\frac{j_3}{Nd_3}\right)}}\;.\nonumber
\end{eqnarray}
Since the indices $j_1,j_2,j_3$ are runing in the range $[-r,r]$ and $1\ll r
\ll N$ then 
\begin{eqnarray}
1+\frac{j_k}{Nd_k}\geq1-\frac{r}{Nd_k}\;,\;\;\;k=1,2,3\;,\nonumber
\end{eqnarray}
that leads to the following inequalities
\begin{eqnarray}
\frac{\Sigma_1(r,N)}{\Sigma_0(r,N)}\leq\prod_{k=1}^3\left(1-\frac{r}{Nd_k}
\right)^{-1/2}\;,\;\;\;\;\frac{\Sigma_2(r,N)}{\Sigma_0(r,N)}\leq\frac{3r-1}{
N(d_1+d_2+d_3)}\prod_{k=1}^3\left(1-\frac{r}{Nd_k}\right)^{-1/2}\;.\label{s4c}
\end{eqnarray}
Combining (\ref{s3c3}) and (\ref{s4c}) we obtain
\begin{eqnarray}
K_{N,r}\left({\bf d}^3\right)>\sqrt{3}-\frac{1}{\sqrt{N}}\left(d_1+d_2+d_3+
\frac{3r-1}{N}\right)\prod_{k=1}^3\left(d_k-\frac{r}{N}\right)^{-1/2}\;,
\nonumber
\end{eqnarray}
and finally the limit yields
\begin{eqnarray}
{\sf K}\left({\bf d}^3\right)\geq\sqrt{3}\;.\label{s5c}
\end{eqnarray}
Thus, Conjecture \ref{con10} is refuted for $m=3$.

As for higher dimension, Conjecture \ref{con10} for numerical semigroups ${\sf 
S}\left({\bf d}^m\right)$, $m\geq 4$, doesn't contradict the best lower bound 
for $C({\bf d}^m)$ known today \cite{kill00}
\begin{eqnarray}  
C\left({\bf d}^m\right)\geq g_m\sqrt[m-1]{d_1\cdot\ldots\cdot d_m}-(d_1+
\ldots+d_m)+1\;.\label{syl2b}
\end{eqnarray}
This leaves Conjecture \ref{con10} open for the case $m\geq 4$.
\subsection{Conjecture $\#$1999--9 and its discussion}\label{ccon2}
Conjecture $\#$1999--9 deals with the asymptotic behavior of the average
distribution of the numerical semigroup ${\sf S}\left({\bf d}^m\right)$ in the 
interval of integers between 0 and $C\left({\bf d}^m\right)$. Denote by $p\left(
{\bf d}^m\right)$ a fraction of the segment {\rm [0;} $C\left({\bf d}^m\right)
-1${\rm ]} which is occupied by the semigroup ${\sf S}\left({\bf d}^m\right)$,
\begin{eqnarray}
p\left({\bf d}^m\right)=\frac{\widetilde{G}\left({\bf d}^m\right)}{C\left(
{\bf d}^m\right)}\;.\label{cal1}
\end{eqnarray}
According to (\ref{intro5}) the fraction $p\left({\bf d}^3\right)$ satisfies 
\begin{eqnarray}
p\left({\bf d}^3\right)&=&\frac{1}{2}\;,\;\mbox{iff}\;\;{\sf S}\left({\bf d}^3
\right)\;\;\mbox{is symmetric}\;,\label{cal0a}\\
p\left({\bf d}^3\right)&<&\frac{1}{2}\;,\;\mbox{iff}\;\;{\sf S}\left({\bf d}^3
\right)\;\;\mbox{is non--symmetric}\;.\label{cal0b}
\end{eqnarray}
\begin{conjecture}{\rm ($\#$1999--9)}\label{con20}
Determine $p\left({\bf d}^m\right)$ for large vectors ${\bf d}^m$. 
Conjecturally, this fraction is asymptotically equal to $1/m$ (with 
overwhelming probability for large ${\bf d}^m$),
\begin{eqnarray}
\widetilde{\sf G}\left({\bf d}^m\right)\stackrel{asymptotically\atop 
weak}{\equiv}\frac{1}{m}{\sf C}\left({\bf d}^m\right)\;.\label{cal1a}
\end{eqnarray}\end{conjecture}
The words ' {\em asymptotically equal} ' and ' {\em with overwhelming
probability for large ${\bf d}^m$} ' presume a weak asymptotics for $p\left(
{\bf d}^m\right)$ via the averaging procedure decribed in Section \ref{ccon1}.

Represent Conjecture \ref{con20} in the sense of (\ref{set1b}),
\begin{eqnarray}
{\sf P}\left({\bf d}^m\right)=\lim_{r,N\to\infty\atop r(N)/N\to 0}p_{N,r}
\left({\bf d}^m\right)=\frac{1}{m}\;,\;\;\;\;\mbox{where}\label{cal2}
\end{eqnarray}
\begin{eqnarray}
p_{N,r}\left({\bf d}^m\right)=\frac{\sum_{j_1,\ldots,j_m=-r}^{r}\widetilde{G}
\left(N{\bf d}^m+{\bf j}^m\right)}{\sum_{j_1,\ldots,j_m=-r}^{r}C\left(N{\bf d}^
m+{\bf j}^m\right)}\;,\;\;\;\mbox{and}\;\;\;\;N{\bf d}^m+{\bf j}^m\in
{\mathbb M}_{N,r}\left({\bf d}^m\right)\;.\label{cal3}
\end{eqnarray}
For $m=2$ the corresponding semigroups ${\sf S}\left(N{\bf d}^2+{\bf j}^2
\right)$ are symmetric and by (\ref{intro5}) we have,
\begin{eqnarray}
{\sf P}\left({\bf d}^2\right)=p_{N,r}\left({\bf d}^2\right)=p\left(N{\bf 
d}^2+{\bf j}^2\right)=p\left({\bf d}^2\right)=\frac{1}{2}\;.\label{den2}
\end{eqnarray}
Consider the 3--dim version of Conjecture \ref{con20} and recall two important 
results which are worthwhile to discuss Conjecture. First, consider a semigroup 
${\sf S}\left({\bf d}^3\right)$ and represent $p_{N,r}\left({\bf d}^3\right)$ 
as follows,
\begin{eqnarray}
p_{N,r}\left({\bf d}^3\right)=\frac{\sum_{j_1,j_2,j_3=-r}^{r}p\left(N{
\bf d}^3+{\bf j}^3\right)C\left(N{\bf d}^3+{\bf j}^3\right)}{\sum_{j_1,j_2,
j_3=-r}^{r}C\left(N{\bf d}^3+{\bf j}^3\right)}\;.\label{den7}
\end{eqnarray}
By Lemma \ref{lem2} a summation in (\ref{den7}) is performed over $N{\bf d}^3+
{\bf j}^3\in {\mathbb M}_{N,r}\left({\bf d}^3\right)$ and the corresponding 
semigroups ${\sf S}\left(N{\bf d}^3+{\bf j}^3\right)$ are non--symmetric only.
Theorems \ref{theo1} and \ref{theo2} imply for such admissible vectors the 
following,
\begin{eqnarray}
p\left(N{\bf d}^3+{\bf j}^3\right)&=&\frac{1}{3}\;,\;\mbox{iff}\;\;N{\bf d}^3
+{\bf j}^3=\{3,3k+1,3k+2\}\;,\;\;k\geq 1\;,\label{den8b}\\
p\left(N{\bf d}^3+{\bf j}^3\right)&>&\frac{1}{3}\;,\;\mbox{otherwise}\;.
\label{den8c}
\end{eqnarray}
However (\ref{den8b}) doesn't hold for any $N$, $d_1$ and $j_1$ due to 
(\ref{s00c}). Thus, the semigroups ${\sf S}\left(N{\bf d}^3+{\bf j}^3\right)$ 
contributing to (\ref{den7}) satisfy (\ref{cal0a}) and (\ref{den8c}), and 
therefore
\begin{eqnarray}
\frac{1}{3}<p_{N,r}\left({\bf d}^3\right)<\frac{1}{2}\;.\label{den9}
\end{eqnarray}
Taking the limit $r,N\to\infty$, $r(N)/N\to 0$ in (\ref{den9}) we refute 
Conjecture \ref{con20} for $m=3$,
\begin{eqnarray}
\frac{1}{3}<{\sf P}\left({\bf d}^3\right)<\frac{1}{2}\;.\label{den20}
\end{eqnarray}
In Section \ref{rev1} we improve the left hand side of inequality (\ref{den20}) 
by applying recent results \cite{fel04} in the Frobenius problem for the 
numerical semigroups ${\sf S}\left({\bf d}^3\right)$.

As for higher dimension, $m\geq 4$, the relations between $G\left({\bf d}^m
\right)$ and $\widetilde{G}\left({\bf d}^m\right)$ do exist \cite{frob87}, 
\cite{brown91} and are similar to those given in Theorems \ref{theo1} and 
\ref{theo2},
\begin{eqnarray}
&&G\left({\bf d}^m\right)\leq \widetilde{G}\left({\bf d}^m\right)
t\left({\sf S}\left({\bf d}^m\right)\right)\;,\nonumber\\
&&G\left({\bf d}^m\right)=\left\{\begin{array}{l}\widetilde{G}\left({\bf d}^m
\right)\;,\;\;\;\mbox{iff}\;\;\;{\sf S}\left({\bf d}^m\right)\;\;\mbox{is 
symmetric}\;,\\
\widetilde{G}\left({\bf d}^m\right)t\left({\sf S}\left({\bf d}^m\right)\right)
\;,\;\;\;\mbox{iff}\;\;\;{\bf d}^m=\{m,km+1,\ldots,km+m-1\}\;,\;\;k\geq 1\;.
\end{array}\right.\nonumber
\end{eqnarray}
However, the type $t\left({\sf S}\left({\bf d}^m\right)\right)$ in the case 
$m\geq 4$ doesn't posses such universal properties as in (\ref{intro6}). 
Here there are very mild constraints only, 
\begin{eqnarray}
t\left({\sf S}\left({\bf d}^m\right)\right)=d_1-1,\;\;\mbox{iff}\;\;d_1=m,\;
\cite{sal79}\;\;\;\;\;\mbox{and}\;\;\;\;\;t\left({\sf S}\left({\bf d}^m\right)
\right)<d_1-1,\;\;\mbox{otherwise},\;\cite{herz71}.\nonumber
\end{eqnarray}
These properties are not enough to resolve Conjecture \ref{con20} for the 
case $m\geq 4$ and leave it open meanwhile.
\subsection{Conjecture $\#$1999--10 and its discussion}\label{cconn3}
Conjecture $\#$1999--10 deals with the asymptotic behavior of the average
distribution of the numerical semigroup ${\sf S}\left({\bf d}^m\right)$ in the 
interval of integers between 0 and $C\left({\bf d}^m\right)$. Examples show 
that semigroup ${\sf S}\left({\bf d}^m\right)$ fills the right half of the 
segment {\rm [0;} $C\left({\bf d}^m\right)-1${\rm ]} more dense.
\begin{conjecture}{\rm ($\#$1999--10)}\label{con30}
Find the typical density of filling the segment {\rm [0;} $C\left({\bf d}^m
\right)-1${\rm ]} asymptotically for large ${\bf d}^m$. The conjectured
behavior of the density $p_m(s)$ at a point $s<C\left({\bf d}^m\right)$ is
\begin{eqnarray}
p_m(s)=\left(\frac{s}{{\sf C}\left({\bf d}^m\right)}\right)^{m-1}\;.
\label{den31}
\end{eqnarray}
Such a distribution would immediately imply that the semigroup ${\sf S}\left({
\bf d}^m\right)$ occupies $1/m$--th part of the segment {\rm [0;} $C\left({\bf 
d}^m\right)-1${\rm ]},
\begin{eqnarray}
\int_0^{{\sf C}\left({\bf d}^m\right)}p_m(s)ds=\frac{{\sf C}\left({\bf d}^m
\right)}{m}\;.\label{den32}
\end{eqnarray}
\end{conjecture}
Since Conjecture \ref{con30} is strongly related by the last Formula 
(\ref{den32}) to Conjecture \ref{con20} and the latter is refuted for $m=3$ in 
Section \ref{ccon2} then the conjectured Formula (\ref{den31}) for $p_3(s)$ is 
not valid.
\section{Conjectures $\#$1999--8 and $\#$1999--9 revisited}\label{rev1}
In Section \ref{arr1} we have refuted Conjectures \ref{con10} and \ref{con20} 
in the case $m=3$ implicitly but have not found the explicit expressions for 
${\sf K}\left({\bf d}^3\right)$ and ${\sf P}\left({\bf d}^3\right)$ although 
both Conjectures ask for them. There is another point which makes our solutions 
in Section \ref{arr1} incomplete. This is an unknown universality of these 
solutions. In other words, do ${\sf K}\left({\bf d}^3\right)$ and 
${\sf P}\left({\bf d}^3\right)$ depend on the vector ${\bf d}^3$ where an 
averaging is performed around, or they are given by real numbers that is 
presumed by Arnol'd ? The question remains actual even in the present 
situation when Conjectures are refuted.

Based on recent results \cite{fel04} in the Frobenius problem for the numerical 
semigroups ${\sf S}\left({\bf d}^3\right)$ we give in this Section the explicit 
expressions for ${\sf K}\left({\bf d}^3\right)$ and ${\sf P}\left({\bf d}^3
\right)$ and show that they are not universal. We also improve an inequality 
(\ref{den20}) by enhancing its lower bound. 
Before going to the subject we recall recent results \cite{fel04} in the 
Frobenius problem for the numerical semigroups ${\sf S}\left({\bf d}^3\right)$.
We focus on non--symmetric semigroups since by Lemma \ref{lem2} such semigroups 
contribute to the values of ${\sf K}\left({\bf d}^3\right)$ and ${\sf P}\left(
{\bf d}^3\right)$ only.
\subsection{Matrix $\widehat {\cal R}_3$ of minimal relations, conductor 
$C\left({\bf d}^3\right)$ and genus $G\left({\bf d}^3\right)$}\label{stat0}
Let ${\sf S}\left(d_1,d_2,d_3\right)\subset {\mathbb Z}_+$ be the additive
numerical semigroup finitely generated by a minimal set of positive integers
$d_1<d_2<d_3$ such that $\gcd(d_1,d_2,d_3)=1$. Following Johnson \cite{john60} 
define {\em the minimal relation} for given triple ${\bf d}^3=(d_1,d_2,d_3)$,
\begin{eqnarray}
&&a_{11}d_1=a_{12}d_2+a_{13}d_3\;,\;\;\;a_{22}d_2=a_{21}d_1+a_{23}d_3\;,\;\;
\;a_{33}d_3=a_{31}d_1+a_{32}d_2\;,\;\;\;\;\mbox{where}\;\;\;\;\;
\label{joh1}\\
&&a_{jj}=\min\left\{v_{jj}\;\bracevert\;v_{jj}\geq 2,\;v_{jj}d_j=v_{jk}d_k+
v_{jl}d_l,\;v_{jk},v_{jl}\in {\mathbb Z}_+\cup\{0\}\right\}\;,\label{joh1a}\\
&&\gcd(a_{jj},a_{jk},a_{jl})=1\;,\;\;\;\mbox{and}\;\;\;(j,k,l)=(1,2,3),\;
(2,3,1),\;(3,1,2)\;,\nonumber
\end{eqnarray}
The uniquely defined values of $v_{ij},i\neq j$ which give $a_{ii}$ will be
denoted by $a_{ij},i\neq j$. Represent (\ref{joh1}) as a matrix equation
\begin{eqnarray}
\widehat {\cal R}_3\left(\begin{array}{r}d_1\\d_2\\d_3 \end{array}\right)=
\left(\begin{array}{r}0\\0\\0\end{array}\right)\;,\;\;\;
\widehat {\cal R}_{3}=\left(\begin{array}{rrr}a_{11}&-a_{12}&-a_{13}\\
-a_{21}&a_{22}&-a_{23}\\-a_{31}&-a_{32} & a_{33}\end{array}\right)\;,\;\;\;
\left\{\begin{array}{r}\gcd(a_{11},a_{12},a_{13})=1\\\gcd(a_{21},a_{22},
a_{23})=1\\\gcd(a_{31},a_{32},a_{33})=1\end{array}\right.\;,\label{joh2}
\end{eqnarray}
and establish {\em a standard form} of the matrix $\widehat {\cal R}_{3}$
satisfying (\ref{joh1}) and (\ref{joh1a}).

For the non--symmetric semigroups the matrix $\widehat {\cal R}_{3}$ can be
written as follows \cite{john60}
\begin{eqnarray}  
\widehat {\cal R}_{3}=\left(\begin{array}{ccc}u_1+w_1&-u_2&-w_3\\-w_1&u_2+w_2&
-u_3\\-u_1&-w_2&u_3+w_3\end{array}\right)\;,\;\;\;
\left\{\begin{array}{l}\gcd(u_1,w_2,u_3+w_3)=1\\\gcd(u_2,w_3,u_1+w_1)=1\\
\gcd(u_3,w_1,u_2+w_2)=1\;,\end{array}\right.\label{joh3}
\end{eqnarray}
where $u_i,w_i\in {\mathbb Z}_+\;,\;\;i=1,2,3$. The generators $d_1$, $d_2$ 
and $d_3$ are uniquely defined in the form \cite{john60}
\begin{eqnarray}
d_1=u_2u_3+w_2w_3+u_2w_3\;,\;\;d_2=u_3u_1+w_3w_1+u_3w_1\;,\;\;
d_3=u_1u_2+w_1w_2+u_1w_2\;.\label{joh4}
\end{eqnarray}
The conductor $C\left({\bf d}^3\right)$ and the genus $G\left({\bf d}^3
\right)$ are given by \cite{fel04}
\begin{eqnarray}
C\left({\bf d}^3\right)&=&1+\prod_{i=1}^3(u_i+w_i)-A_2-B_2-(u_1w_2+u_2w_3+
u_3w_1)+max\{A_3,B_3\}\;,\label{joh5}\\
2G\left({\bf d}^3\right)&=&1+\prod_{i=1}^3(u_i+w_i)-A_2-B_2-(u_1w_2+u_2w_3+
u_3w_1)+A_3+B_3\;,\;\;\;\mbox{where}\label{joh6}
\end{eqnarray}
$$
A_2=u_1u_2+u_3u_1+u_2u_3\;,\;\;A_3=u_1u_2u_3\;,\;\;B_2=w_1w_2+w_3w_1+w_2w_3
\;,\;\;B_3=w_1w_2w_3\;.
$$ 
Notice that
\begin{eqnarray}
2G\left({\bf d}^3\right)-C\left({\bf d}^3\right)=min\{A_3,B_3\}\;.\label{joh7}
\end{eqnarray}
\subsection{Explicit expression for ${\sf K}\left({\bf d}^3\right)$ and 
its lower bound}\label{low10}
In \cite{arn04} Arnol'd gave a weak version of Conjecture \ref{con10}:

{\em For growing values of N and $r$, $r(N)/N\to 0$ when $N\to \infty$, and 
large ${\bf d}^m$ the mean values $C_{N,r}\left({\bf d}^m\right)$ have a limit} 
({\em probably provided by conjectured formula} (\ref{syl1})) which grows as 
\begin{eqnarray}
const\sqrt[m-1]{d_1\cdot\ldots\cdot d_m}\;.\label{joh8}
\end{eqnarray}
Here a conjectured limit (\ref{joh8}) is more weak than (\ref{syl1}) since 
it admits $const\neq g_m$. Although it does claim the similar dependence
$\sqrt[m-1]{d_1\cdot\ldots\cdot d_m}$ as Conjecture \ref{con10} does. In that 
sense our solution in Section \ref{ccon1} refutes (\ref{syl1}) but its weak 
version (\ref{joh8}) remains still open.

Consider a non--symmetric semigroup ${\sf S}\left({\bf d}^3\right)$ and
calculate a function $K_{N,r}\left({\bf d}^3\right)$ given in (\ref{s0c}).
Formulas (\ref{joh4}) and (\ref{joh5}) dictate to perform the averaging of 
numerical function $A\left({\bf d}^3\right)$ not on the usual 3--dim cubic 
lattice ${\mathbb Z}_+^3$, where a set ${\mathbb M}_{N,r}\left({\bf d}^3
\right)$ is defined by (\ref{set1}), but on the cubic lattice of higher 
dimension. Namely, denote by ${\bf u}^3$ and ${\bf w}^3$ two 3--dim tuples 
$(u_1,u_2,u_3)$ and $(w_1,w_2,w_3)$, respectively. Consider their union ${\bf 
u}^3\cup {\bf w}^3=(u_1,u_2,u_3,w_1,w_2,w_3)$ as a tuple in the 6--dim cubic 
lattice ${\mathbb Z}_+^3\times {\mathbb Z}_+^3$ as follows,
\begin{eqnarray}
{\mathbb Z}_+^3\times {\mathbb Z}_+^3:=\left\{{\bf u}^3\cup {\bf w}^3\;
\bracevert\;{\bf u}^3\cup {\bf w}^3=(u_1,u_2,u_3,w_1,w_2,w_3)\;,\;u_i,w_i\in
{\mathbb Z}_+\right\}\;.\label{k1}
\end{eqnarray}
A mapping ${\mathbb Z}_+^3\times {\mathbb Z}_+^3\longmapsto {\mathbb Z}_+^3$ 
is defined by equations (\ref{joh4}). In order to find a weak asymptotics 
replace a scaling in ${\mathbb Z}^3$ lattice, $N^2d_i\in {\mathbb Z}_+,\;N\in{
\mathbb Z}_+$, by the scaling in ${\mathbb Z}^3\times {\mathbb Z}_+^3$ lattice, 
$Nu_i,Nw_i\in {\mathbb Z}_+$, and define a set ${\mathbb A}_{N,r}\left({\bf u}^
3\cup {\bf w}^3\right)$ on ${\mathbb Z}_+^3\times {\mathbb Z}_+^3$ as follows,
\begin{eqnarray}
{\mathbb A}_{N,r}\left({\bf u}^3\cup {\bf w}^3\right)=\left\{\left(N{\bf u}^3+
{\bf j}^3\right)\cup \left(N{\bf w}^3+{\bf k}^3\right)\left\bracevert\left.
\begin{array}{c}\gcd\left(D_{1,N}(j_i,k_i),D_{2,N}(j_i,k_i),D_{3,N}(j_i,k_i)
\right)=1\;,\\{\bf j}^3=(j_1,j_2,j_3)\;,\;\;\;{\bf k}^3=(k_1,k_2,k_3)\;,\\-r
\leq j_i,k_i\leq r\;,\;\;1\ll r\ll N\end{array}\right.\right.\right\}\nonumber
\end{eqnarray}
where 
\begin{eqnarray}
D_{1,N}(j_i,k_i)=(Nu_2+j_2)(Nu_3+j_3)+(Nw_2+k_2)(Nw_3+k_3)+(Nu_2+j_2)(Nw_3+k_3)
\;,\nonumber\\
D_{2,N}(j_i,k_i)=(Nu_3+j_3)(Nu_1+j_1)+(Nw_3+k_3)(Nw_1+k_1)+(Nu_3+j_3)(Nw_1+k_1)
\;,\nonumber\\
D_{3,N}(j_i,k_i)=(Nu_1+j_1)(Nu_2+j_2)+(Nw_1+k_1)(Nw_2+k_2)+(Nu_1+j_1)(Nw_2+k_2)
\;.\nonumber
\end{eqnarray}
A cardinality of a set ${\mathbb A}_{N,r}\left({\bf u}^3\cup {\bf w}^3
\right)$ can be estimated in the same way as was done in Section \ref{met1} 
for the set ${\mathbb M}_{N,r}\left({\bf d}^3\right)$, 
\begin{eqnarray}
\#\left\{{\mathbb A}_{N,r}\left({\bf u}^3\cup {\bf w}^3\right)\right\}\simeq
\frac{(2r)^6}{\zeta(3)}\;.\label{k1a}
\end{eqnarray}
Substituting (\ref{joh4}) and (\ref{joh5}) into (\ref{s0c}) and averaging over 
the set ${\mathbb A}_{N,r}\left({\bf u}^3\cup {\bf w}^3\right)$ we get
\begin{eqnarray}
K_{N,r}\left({\bf d}^3\right)=\frac{\sum_{j_1,j_2,j_3\atop k_1,k_2,k_3}^{
{\mathbb A}_{N,r}}C_{j_1,j_2,j_3}^{k_1,k_2,k_3}}{\sum_{j_1,j_2,j_3\atop k_1,
k_2,k_3}^{{\mathbb A}_{N,r}}\sqrt{V_{j_1,j_2,j_3}^{k_1,k_2,k_3}}}\;,\;\;\;\;
\mbox{where}\label{k2}
\end{eqnarray}
\begin{eqnarray}
C_{j_1,j_2,j_3}^{k_1,k_2,k_3}&=&\frac{1}{N^3}+
\left(u_1+w_1+\frac{j_1+k_1}{N}\right)\left(u_2+w_2+\frac{j_2+k_2}{N}\right)
\left(u_3+w_3+\frac{j_3+k_3}{N}\right)+\label{k2a}\\
&&\max\left\{\left(u_1+\frac{j_1}{N}\right)\left(u_2+\frac{j_2}{N}\right)
\left(u_3+\frac{j_3}{N}\right),\left(w_1+\frac{k_1}{N}\right)\left(w_2+
\frac{k_2}{N}\right)\left(w_3+\frac{k_3}{N}\right)\right\}-\nonumber\\
&&\frac{1}{N}\left[\left(u_1+\frac{j_1}{N}\right)\left(u_2+ 
\frac{j_2}{N}\right)+\left(u_3+\frac{j_3}{N}\right)
\left(u_1+\frac{j_1}{N}\right)+\left(u_2+\frac{j_2}{N}\right)
\left(u_3+\frac{j_3}{N}\right)\right]-\nonumber\\
&&\frac{1}{N}\left[\left(w_1+\frac{k_1}{N}\right)\left(w_2+\frac{k_2}{N}\right)+
\left(w_3+\frac{k_3}{N}\right)\left(w_1+\frac{k_1}{N}\right)+\left(w_2+
\frac{k_2}{N}\right)\left(w_3+\frac{k_3}{N}\right)\right]-\nonumber\\
&&\frac{1}{N}\left[\left(u_1+\frac{j_1}{N}\right)\left(w_2+\frac{k_2}{N}\right)+
\left(u_2+\frac{j_2}{N}\right)\left(w_3+\frac{k_3}{N}\right)+\left(u_3+
\frac{j_3}{N}\right)\left(w_1+\frac{k_1}{N}\right)\right]\;,\nonumber
\end{eqnarray}
\begin{eqnarray}
V_{j_1,j_2,j_3}^{k_1,k_2,k_3}&=&\left[\left(u_2+\frac{j_2}{N}\right)
\left(u_3+\frac{j_3}{N}\right)+\left(w_2+\frac{k_2}{N}\right)
\left(w_3+\frac{k_3}{N}\right)+\left(u_2+\frac{j_2}{N}\right)
\left(w_3+\frac{k_3}{N}\right)\right]\times\nonumber\\
&&\left[\left(u_3+\frac{j_3}{N}\right)\left(u_1+\frac{j_1}{N}\right)+
\left(w_3+\frac{k_3}{N}\right)\left(w_1+\frac{k_1}{N}\right)+
\left(u_3+\frac{j_3}{N}\right)\left(w_1+\frac{k_1}{N}\right)\right]
\times\nonumber\\
&&\left[\left(u_1+\frac{j_1}{N}\right)\left(u_2+\frac{j_2}{N}\right)+
\left(w_1+\frac{k_1}{N}\right)\left(w_2+\frac{k_2}{N}\right)+
\left(u_1+\frac{j_1}{N}\right)\left(w_2+\frac{k_2}{N}\right)\right]\;.\nonumber
\end{eqnarray}
An upper limit in (\ref{k2}) means that a summation is performed for $\left(N{
\bf u}^3+{\bf j}^3\right)\cup \left(N{\bf w}^3+{\bf k}^3\right)\in {\mathbb A}
_{N,r}\left({\bf u}^3\cup {\bf w}^3\right)$. Bearing in mind that 
$\sum_{j_1,j_2,j_3\atop k_1,k_2,k_3}^{{\mathbb A}_{N,r}}1=\#\left\{{\mathbb A}_
{N,r}\left({\bf u}^3\cup {\bf w}^3\right)\right\}$ and estimating the terms,
\begin{eqnarray}
\left\vert\left.\sum_{j_1,j_2,j_3\atop k_1,k_2,k_3}^{{\mathbb A}_{N,r}}
\frac{j_i}{N}\right.\right\vert\simeq
\left\vert\left.\sum_{j_1,j_2,j_3\atop k_1,k_2,k_3}^{{\mathbb A}_{N,r}}
\frac{k_i}{N}\right.\right\vert<\frac{2r}{N}(2r)^6,\;\;
\left\vert\left.\sum_{j_1,j_2,j_3\atop k_1,k_2,k_3}^{{\mathbb A}_
{N,r}}\frac{j_ij_l}{N^2}\right.\right\vert\simeq
\left\vert\left.\sum_{j_1,j_2,j_3\atop k_1,k_2,k_3}^{{\mathbb A}_
{N,r}}\frac{j_ik_l}{N^2}\right.\right\vert\simeq
\left\vert\left.\sum_{j_1,j_2,j_3\atop k_1,k_2,k_3}^{{\mathbb A}_
{N,r}}\frac{k_ik_l}{N^2}\right.\right\vert
<\frac{(2r)^2}{N^2}(2r)^6,\nonumber\\
\left\vert\left.\sum_{j_1,j_2,j_3\atop k_1,k_2,k_3}^{{\mathbb A}_{N,r}}
\frac{j_ij_lj_n}{N^3}\right.\right\vert\simeq
\left\vert\left.\sum_{j_1,j_2,j_3\atop k_1,k_2,k_3}^{{\mathbb A}_{N,r}}
\frac{j_ij_lk_n}{N^3}\right.\right\vert\simeq
\left\vert\left.\sum_{j_1,j_2,j_3\atop k_1,k_2,k_3}^{{\mathbb A}_{N,r}}
\frac{j_ik_lk_n}{N^3}\right.\right\vert\simeq
\left\vert\left.\sum_{j_1,j_2,j_3\atop k_1,k_2,k_3}^{{\mathbb A}_{N,r}}
\frac{k_ik_lk_n}{N^3}\right.\right\vert<\frac{(2r)^3}{N^3}(2r)^6\;,\;\;\;
\mbox{etc},\nonumber
\end{eqnarray}
we arrive in accordance with (\ref{k1a}) to the leading terms ${\sf K_n}(u_i,
w_i)$ and ${\sf K_d}(u_i,w_i)$ which are contributing to the both sums in 
(\ref{k2}) in the limit $r,N\to\infty$, $r(N)/N\to 0$,
\begin{eqnarray}
\frac{\sum_{j_1,j_2,j_3\atop k_1,k_2,k_3}^{{\mathbb A}_{N,r}}C_{j_1,j_2,j_3}^
{k_1,k_2,k_3}}{\#\left\{{\mathbb A}_{N,r}\left({\bf u}^3\cup {\bf w}^3\right)
\right\}}\simeq {\sf K_n}(u_i,w_i)+{\cal O}\left(\frac{r}{N}\right),\;\;\;
\frac{\sum_{j_1,j_2,j_3\atop k_1,k_2,k_3}^{{\mathbb A}_{N,r}}\sqrt{V_{j_1,j_2,
j_3}^{k_1,k_2,k_3}}}{\#\left\{{\mathbb A}_{N,r}\left({\bf u}^3\cup {\bf w}^3
\right)\right\}}\simeq {\sf K_d}(u_i,w_i)+{\cal O}\left(\frac{r}{N}\right),
\label{k2b}
\end{eqnarray}
where 
\begin{eqnarray}
{\sf K_n}(u_i,w_i)&=&(u_1+w_1)(u_2+w_2)(u_3+w_3)+max\{u_1u_2u_3,w_1w_2w_3\}\;,
\nonumber\\
{\sf K_d}(u_i,w_i)&=&\sqrt{(u_2u_3+w_2w_3+u_2w_3)(u_3u_1+w_3w_1+u_3w_1)
(u_1u_2+w_1w_2+u_1w_2)}\;.\nonumber
\end{eqnarray}

Finally, we obtain the expression for ${\sf K}\left({\bf d}^3\right)$ in 
accordance with (\ref{s1c})
\begin{eqnarray}
{\sf K}\left({\bf d}^3\right)=\frac{(u_1+w_1)(u_2+w_2)(u_3+w_3)+max\{u_1u_2u_3,
w_1w_2w_3\}}{\sqrt{(u_2u_3+w_2w_3+u_2w_3)(u_3u_1+w_3w_1+u_3w_1)(u_1u_2+w_1w_2+
u_1w_2)}}\;.\label{k3}
\end{eqnarray}
One can show that ${\sf K}\left({\bf d}^3\right)$ attains its minimal value, 
${\sf K}\left({\bf d}^3\right)=\sqrt{3}$ when $u_1=w_1$, $u_2=w_2$ and $u_3=w_3$
(see Appendix \ref{appendix1}). This nicely concides with (\ref{s5c}). The 
representations (\ref{k3}) tells one more important thing: ${\sf K}\left({\bf 
d}^3\right)$ is not universal and depends on the vector ${\bf d}^3$ where an 
averaging is performed around. This refutes Conjecture \ref{con10} in its weak 
version (\ref{joh8}).
\subsection{Explicit expression for ${\sf P}\left({\bf d}^3\right)$ and its 
lower bound}\label{low1}
Consider the non--symmetric semigroup ${\sf S}\left({\bf d}^3\right)$ and 
define a new ratio,
\begin{eqnarray}
q\left({\bf d}^3\right):=\frac{G\left({\bf d}^3\right)-\widetilde{G}\left(
{\bf d}^3\right)}{C\left({\bf d}^3\right)}\;.\label{syl4}
\end{eqnarray}
Associate with it a corresponding function $q_{N,r}\left({\bf d}^m\right)$,
\begin{eqnarray} 
q_{N,r}\left({\bf d}^m\right)=\frac{\sum_{j_1,\ldots,j_m=-r}^{r}\left(G\left(
N{\bf d}^m+{\bf j}^m\right)-\widetilde{G}\left(N{\bf d}^m+{\bf j}^m\right)
\right)}{\sum_{j_1,\ldots,j_m=-r}^{r}C\left(N{\bf d}^m+{\bf j}^m\right)}\;,\;
\;\;\mbox{where}\;\;\;\;N{\bf d}^m+{\bf j}^m\in {\mathbb M}_{N,r}\left({\bf 
d}^m\right).\nonumber
\end{eqnarray}
Both fractions, $p\left({\bf d}^3\right)$ and $q\left({\bf d}^3\right)$, are 
readily related to each other,
\begin{eqnarray}
p\left({\bf d}^3\right)=\frac{1}{2}\left(1-q\left({\bf d}^3\right)\right)\;.
\label{syl4u}
\end{eqnarray}
The same relation holds for theier weak asymptotics,
\begin{eqnarray}
{\sf P}\left({\bf d}^3\right)=\frac{1}{2}\left(1-{\sf Q}\left({\bf d}^3\right)
\right)\;,\;\;\;\mbox{where}\;\;\;{\sf Q}\left({\bf d}^3\right)=\lim_{r,N\to
\infty\atop r(N)/N\to 0}q_{N,r}\left({\bf d}^m\right)\;.\label{syl40}
\end{eqnarray}
Perform the averaging of $q_{N,r}\left({\bf d}^3\right)$ over the set ${\mathbb 
A}_{N,r}\left({\bf u}^3\cup {\bf w}^3\right)$ on the 6--dim cubic lattice 
${\mathbb Z}_+^3\times {\mathbb Z}_+^3$ in terms of the $\widehat {\cal R}_3$ 
matrix entries in the same way as was done in Section \ref{low10} for $K_{N,r}
\left({\bf d}^3\right)$. Bearing in mind (\ref{joh7}) we have,
\begin{eqnarray}
q_{N,r}\left({\bf d}^m\right)=\frac{\sum_{j_1,j_2,j_3\atop k_1,k_2,k_3}
^{{\mathbb A}_{N,r}}M_{j_1,j_2,j_3}^{k_1,k_2,k_3}}{\sum_{j_1,j_2,j_3\atop 
k_1,k_2,k_3}^{{\mathbb A}_{N,r}}C_{j_1,j_2,j_3}^{k_1,k_2,k_3}}\;,\;\;\;\;\mbox{where}\;\;\;
\left(N{\bf u}^3+{\bf j}^3\right)\cup \left(N{\bf w}^3+{\bf k}^3\right)\in
{\mathbb A}_{N,r}\left({\bf u}^3\cup {\bf w}^3\right).\label{syl80}
\end{eqnarray}
A denominator of (\ref{syl80}) is defined in (\ref{k2a}), and a numerator reads,
\begin{eqnarray}
M_{j_1,j_2,j_3}^{k_1,k_2,k_3}=\min\left\{\left(u_1+\frac{j_1}{N}\right)
\left(u_2+\frac{j_2}{N}\right)\left(u_3+\frac{j_3}{N}\right),
\left(w_1+\frac{k_1}{N}\right)\left(w_2+\frac{k_2}{N}\right)
\left(w_3+\frac{k_3}{N}\right)\right\}\;.\nonumber
\end{eqnarray}
Performing summation in (\ref{syl80}) by applying the similar considerations 
as in Section \ref{low10} and taking the limit $r,N\to\infty$, $r(N)/N\to 0$ 
we get finally,
\begin{eqnarray}
{\sf Q}\left({\bf d}^3\right)=\frac{min\{u_1u_2u_3,w_1w_2w_3\}}{(u_1+w_1)
(u_2+w_2)(u_3+w_3)+\max\{u_1u_2u_3,w_1w_2w_3\}}\;.\label{syl90}
\end{eqnarray}
Represent (\ref{syl90}) as follows,
\begin{eqnarray}
{\sf Q}\left({\bf d}^3\right)=\frac{min\{1,\rho_1\rho_2\rho_3\}}{(1+\rho_1)(1+
\rho_2)(1+\rho_3)+max\{1,\rho_1\rho_2\rho_3\}}\;,\;\;\;\rho_i=\frac{u_i}{w_i}
\;,\;\;\;0<\rho_i<\infty\;.\nonumber 
\end{eqnarray}
Making use of inequalities \cite{hardy59} for three basic polynomial invariants 
$\Gamma_1=\rho_1+\rho_2+\rho_3$, $\Gamma_2=\rho_1\rho_2+\rho_2\rho_3+\rho_3
\rho_1$ and $\Gamma_3=\rho_1\rho_2\rho_3$ of symmetric group $S_3$ acting on 
the set $\{\rho_1,\rho_2,\rho_3\}$, 
\begin{eqnarray}
\Gamma_1\geq\sqrt{3\Gamma_2}\geq3\sqrt[3]{\Gamma_3}\;,\nonumber
\end{eqnarray}
we get ${\sf Q}\left({\bf d}^3\right)>0$ and 
\begin{eqnarray}
{\sf Q}\left({\bf d}^3\right)&=&\frac{1}{1+\Gamma_1+\Gamma_2+2\Gamma_3}<
\frac{1}{9}\;,\;\;\;\mbox{if}\;\;\;\Gamma_3\geq 1\;,\label{syl92a}\\
{\sf Q}\left({\bf d}^3\right)&=&\frac{\Gamma_3}{2+\Gamma_1+\Gamma_2+
\Gamma_3}=\frac{1}{1+\Gamma_1\Gamma_3^{-1}+\Gamma_2\Gamma_3^{-1}+2\Gamma_3^
{-1}}<\frac{1}{9}\;,\;\;\;\mbox{if}\;\;\;\Gamma_3\leq 1\;,\label{syl92b}
\end{eqnarray}
since $\Gamma_2\Gamma_3^{-1}=\rho_1^{-1}+\rho_2^{-1}+\rho_3^{-1}$, $\Gamma_1
\Gamma_3^{-1}=\rho_1^{-1}\rho_2^{-1}+\rho_2^{-1}\rho_3^{-1}+\rho_3^{-1}\rho_1^
{-1}$ and $\Gamma_3^{-1}=\rho_1^{-1}\rho_2^{-1}\rho_3^{-1}$ can be considered 
as basic polynomial invariants of symmetric group $S_3$ acting on the set 
$\{\rho_1^{-1},\rho_2^{-1},\rho_3^{-1}\}$.

The case $\rho_1=\rho_2=\rho_3=1$ has to be excluded since the corresponding 
matrix of minimal relations has the entries $u_1=w_1$, $u_2=w_2$ and $u_3=w_3$ 
that results in $\gcd(d_1,d_2,d_3)=3$. This is why both inequalities in 
(\ref{syl92a}) and (\ref{syl92b}) are rigorous. Finally, we obtain by 
(\ref{syl40}) the lower and upper bounds for ${\sf P}\left({\bf d}^3\right)$,
\begin{eqnarray}
\frac{4}{9}<{\sf P}\left({\bf d}^3\right)<\frac{1}{2}\;.\label{syl93}
\end{eqnarray}
The representations (\ref{syl90}) claims that ${\sf P}\left({\bf d}^3\right)$ 
is not universal and depends on the vector ${\bf d}^3$ where an averaging is 
performed around.
\section*{Acknowledgement}
The usefull discussions with A. Juhasz are highly appreciated. 
\appendix
\renewcommand{\theequation}{\thesection\arabic{equation}}
\section{Lower bound of ${\sf K}\left({\bf d}^3\right)$}
\label{appendix1}
\setcounter{equation}{0}
Represent the function ${\sf K}\left({\bf d}^3\right)$ given in (\ref{k3}) as 
follows,
\begin{eqnarray}
{\sf K}\left({\bf d}^3\right)=\frac{(1+\rho_1)(1+\rho_2)(1+\rho_3)+
max\{1,\rho_1\rho_2\rho_3\}}{\sqrt{(1+\rho_2\rho_3+\rho_2)(1+\rho_3\rho_1+
\rho_3)(1+\rho_1\rho_2+\rho_1)}}\;,\;\;\;\rho_i=\frac{u_i}{w_i}\;,\;\;\;0<
\rho_i<\infty\;,\label{app7}
\end{eqnarray}
and consider its square, ${\sf K}^2\left({\bf d}^3\right)={\sf L}\left(\rho_1,
\rho_2,\rho_3\right)$.

First, prove that ${\sf L}\left(\rho_1,\rho_2,\rho_3\right)$ is unbounded 
from above. Consider $0<\overline{\rho_i}<\infty$ such that $\overline{\rho_1}
=1/(\overline{\rho_2}\;\overline{\rho_3})$ and $\overline{\rho_2},\overline{
\rho_3}\gg 1,\overline{\rho_1}\ll 1$. Calculate a leading term in ${\sf L}
\left(\overline{\rho_1},\overline{\rho_2},\overline{\rho_3}\right)$
\begin{eqnarray}
{\sf L}\left(\overline{\rho_1},\overline{\rho_2},\overline{\rho_3}\right)
=\frac{\left[(1+\overline{\rho_1})(1+\overline{\rho_2})(1+\overline{\rho_3})
+1\}\right]^2}{(1+\overline{\rho_2}\;\overline{\rho_3}+\overline{\rho_2})(1+
\overline{\rho_3}\;\overline{\rho_1}+\overline{\rho_3})(1+\overline{\rho_1}
\;\overline{\rho_2}+\overline{\rho_1})}\simeq\frac{\overline{\rho_2}^2
\overline{\rho_3}^2}{\overline{\rho_2}^2\overline{\rho_3}^2\overline{
\rho_1}}=\frac{1}{\overline{\rho_1}}\gg 1\;.\label{app7a}
\end{eqnarray}
The last inequality proves a statement.

Observe that ${\sf L}\left(\rho_1,\rho_2,\rho_3\right)$ is invariant under 
cyclic permutation of variables $\rho_1,\rho_2,\rho_3$,
\begin{eqnarray}  
{\sf L}\left(\rho_1,\rho_2,\rho_3\right)={\sf L}\left(\rho_2,\rho_3,\rho_1
\right)={\sf L}\left(\rho_3,\rho_1,\rho_2\right)\;,\label{app7b}
\end{eqnarray}
and can be represented in four polynomial invariants $\Gamma_i$ of the cyclic 
group $C_3$ \cite{Benson93},
\begin{eqnarray}
&&\Gamma_1=\rho_1+\rho_2+\rho_3\;,\;\;\Gamma_2=\rho_1\rho_2+\rho_2\rho_3+
\rho_3\rho_1\;,\;\;\Gamma_3=\rho_1\rho_2\rho_3\;,\label{app8a}\\
&&\Gamma_4=(\rho_1-\rho_2)(\rho_2-\rho_3)(\rho_3-\rho_1)\;,\;\;\;\mbox{where}
\;\;\Gamma_4^2=\Gamma_1^2\Gamma_2^2+18\Gamma_1\Gamma_2\Gamma_3-4\Gamma_2^3-
4\Gamma_1^3\Gamma_3-27\Gamma_3^2\;.\nonumber
\end{eqnarray}
In both regions, $\Gamma_3> 1$ and $\Gamma_3< 1$, the function ${\sf L}\left(
\rho_1,\rho_2,\rho_3\right)$ is differentiable and attains its extremal values 
if $\partial {\sf L}/\partial\rho_i=0,\;i=1,2,3$. In other words, 
\begin{eqnarray}
\frac{\partial {\sf L}}{\partial\rho_i}=\sum_{j=1}^4\frac{\partial {\sf L}}
{\partial\Gamma_j}\frac{\partial \Gamma_j}{\partial\rho_i}=0\;,\;\;\;i=1,2,3
\;,\nonumber
\end{eqnarray}
or
\begin{eqnarray}
\left(\frac{\partial {\sf L}}{\partial\Gamma_1}+\frac{\partial {\sf L}}
{\partial\Gamma_4}\frac{\partial\Gamma_4}{\partial\Gamma_1}\right)
\frac{\partial\Gamma_1}{\partial\rho_i}+
\left(\frac{\partial {\sf L}}{\partial\Gamma_2}+\frac{\partial {\sf L}}
{\partial\Gamma_4}\frac{\partial\Gamma_4}{\partial\Gamma_2}\right)
\frac{\partial\Gamma_2}{\partial\rho_i}+
\left(\frac{\partial {\sf L}}{\partial\Gamma_3}+\frac{\partial {\sf L}}
{\partial\Gamma_4}\frac{\partial\Gamma_4}{\partial\Gamma_3}\right)
\frac{\partial\Gamma_3}{\partial\rho_i}=0\;.\label{app9}
\end{eqnarray}
Substituting (\ref{app8a}) into (\ref{app9}) and removing singular 
multiplier $1/\Gamma_4$ after taking derivatives $\partial\Gamma_4/\partial
\Gamma_i$ we get three equations for $i=1,2,3$,
\begin{eqnarray}
\left(\Gamma_4\frac{\partial {\sf L}}{\partial\Gamma_1}+K_{41}\frac{
\partial {\sf L}}{\partial\Gamma_4}\right)\frac{\partial\Gamma_1}{\partial
\rho_i}+\left(\Gamma_4\frac{\partial {\sf L}}{\partial\Gamma_2}+K_{42}
\frac{\partial {\sf L}}{\partial\Gamma_4}\right)\frac{\partial\Gamma_2}{
\partial\rho_i}+\left(\Gamma_4\frac{\partial {\sf L}}{\partial\Gamma_3}+
K_{43}\frac{\partial {\sf L}}{\partial\Gamma_4}\right)\frac{\partial
\Gamma_3}{\partial\rho_i}=0\;,\label{app10}
\end{eqnarray}
where
\begin{eqnarray}
K_{41}=\Gamma_1\Gamma_2^2+9\Gamma_2\Gamma_3-6\Gamma_1^2\Gamma_3\;,\;\;\;
K_{42}=\Gamma_1^2\Gamma_2+9\Gamma_1\Gamma_3-6\Gamma_2^2\;,\;\;\;
K_{43}=9\Gamma_1\Gamma_2-2\Gamma_1^3-27\Gamma_3\;.\nonumber
\end{eqnarray}
Equations (\ref{app10}) have nontrivial solutions if $\det\left(\left(
\partial\Gamma_j/\partial\rho_i\right)\right)=0$. Substituting (\ref{app8a})
into the last equality we obtain $\det\left(\left(\partial\Gamma_j/\partial
\rho_i\right)\right)=\Gamma_4=0$. In other words, ${\sf L}\left(\rho_1,\rho_2,\rho_3
\right)$ attains its extremum at the planes $\rho_1=\rho_2$, $\rho_2=
\rho_3$ and $\rho_3=\rho_1$ where 
$$
K_{4i}\left(\rho_1,\rho_2,\rho_2\right)\neq 0\;,\;\;\;K_{4i}\left(\rho_3,
\rho_2,\rho_3\right)\neq 0\;,\;\;\;K_{4i}\left(\rho_1,\rho_1,\rho_3\right)
\neq 0\;,\;\;\;\;\;i=1,2,3\;.
$$ 
A cyclic invariance (\ref{app7b}) of ${\sf L}\left(\rho_1,\rho_2,\rho_3\right)$ 
makes all three planes equivalent in the sense that provides for ${\sf L}\left(
\rho_1,\rho_2,\rho_3\right)$ the same kind of extremum which can be only a 
minimum due to (\ref{app7a}). Consider one of the solutions, when $\rho_1\neq
\rho_2=\rho_3$, 
\begin{eqnarray}
{\sf L}\left(\rho_1,\rho_2,\rho_2\right)=\frac{\left[(1+\rho_1)(1+
\rho_2)^2+max\{1,\rho_1\rho_2^2\}\right]^2}{(1+\rho_1\rho_2+\rho_1)(1+
\rho_1\rho_2+\rho_2)(1+\rho_2^2+\rho_2)}\;.\label{app11}
\end{eqnarray}
This function possesses additional invariance under inversion of both variables,
\begin{eqnarray}
{\sf L}\left(\rho_1,\rho_2,\rho_2\right)={\sf L}\left(\frac{1}{\rho_1},
\frac{1}{\rho_2},\frac{1}{\rho_2}\right)\;.\label{app12}
\end{eqnarray}
Last relation (\ref{app12}) simplifies essentialy further consideration since 
if in a region $\rho_1\rho_2^2\geq 1$ holds inequality ${\sf L}\left(\rho_1,
\rho_2,\rho_2\right)>const$ then it also holds in a region $\rho_1\rho_2^2\leq
1$.

Consider a region $\rho_1\rho_2^2\geq 1$ and represent ${\sf L}\left(\rho_1,
\rho_2,\rho_2\right)$ as follows,
\begin{eqnarray}
{\sf L}\left(\rho_1,\rho_2,\rho_2\right)=3+\frac{{\cal U}^2-3{\cal W}}{{\cal W}}
\;,\;\;\;\;\mbox{where}\label{app13}
\end{eqnarray}
\begin{eqnarray}
{\cal U}=(1+\rho_1)(1+\rho_2)^2+\rho_1\rho_2^2\;,\;\;\;\;{\cal W}=(1+\rho_1
\rho_2+\rho_1)(1+\rho_1\rho_2+\rho_2)(1+\rho_2^2+\rho_2)\;.\nonumber
\end{eqnarray}
Denoting $\rho_1=\varepsilon+1/\rho_2^2,\;\varepsilon\geq 0$, we get
\begin{eqnarray}
\frac{{\cal U}^2-3{\cal W}}{{\cal W}}=\frac{1}{\rho_2}\;\frac{{\cal T}_0+
{\cal T}_1\varepsilon+{\cal T}_2\varepsilon^2}{{\cal P}_0+{\cal P}_1
\varepsilon+{\cal P}_2\varepsilon^2}\;,\;\;\;\;\;\;\mbox{where}\label{app14}
\end{eqnarray}
\begin{eqnarray}
&&{\cal T}_0=(\rho_2-1)^2(1+\rho_2+\rho_2^2)^3\;,\;\;\;{\cal T}_1=\rho_2^2
\left(1+(\rho_2-1)^2\right)(1+\rho_2+\rho_2^2)^2\;,\nonumber\\
&&{\cal T}_2=\rho_2^4\left(1+\rho_2(1+\rho_2)(1+\rho_2+\rho_2^2)\right)\;,
\label{app15}\\
&&{\cal P}_0=1+3\rho_2+6\rho_2^2+7\rho_2^3+6\rho_2^4+3\rho_2^5+\rho_2^6\;,\;
\;\;{\cal P}_1=\rho_2^2\left(2+5\rho_2+8\rho_2^2+7\rho_2^3+4\rho_2^4+\rho_2^5
\right)\;,\nonumber\\
&&{\cal P}_2=\rho_2^4\left(1+2\rho_2+2\rho_2^2+\rho_2^3\right)\;.
\label{app16}
\end{eqnarray}
One can see from (\ref{app15}) and (\ref{app16}) that ${\cal P}_0\geq 1,
{\cal P}_1,{\cal P}_2,{\cal T}_0,{\cal T}_1,{\cal T}_2\geq 0$. Thus, 
${\sf L}\left(\rho_1,\rho_2,\rho_2\right)$ arrives its minimal value 3 
when $\varepsilon=0$, $\rho_2=1$. So the whole function ${\sf L}\left(
\rho_1,\rho_2,\rho_3\right)$ also arrives its minimal value 3 when 
$\rho_1=\rho_2=\rho_3=1$.

\end{document}